\documentclass[11pt]{article}
\usepackage{psfig,amssymb,amsmath}
\usepackage{epic}
\newtheorem{theorem}{Theorem}[section]
\newtheorem{cor}[theorem]{Corollary}
\newtheorem{defi}{Definition}[section]
\newtheorem{lemma}[theorem]{Lemma}
\newtheorem{prop}[theorem]{Proposition}

\newcommand{\bsq}{\vrule height .9ex width .8ex depth -.1ex}

\newcommand{\hsp}{\hspace*{\parindent}}

\newcommand{\sgn}{{\mbox{sgn}}}

\newcommand{\af}{\alpha}
\newcommand{\ep}{\epsilon}

\newcommand{\vu}{\nu}

\newcommand{\RR}{{\mathbb R}}
\newcommand{\QQ}{{\mathbb Q}}
\newcommand{\TT}{{\bf T}}
\newcommand{\ZZ}{{\mathbb Z}}
\newcommand{\CC}{{\mathbb C}}

\newcommand{\sga}{{\mathfrak a}}
\newcommand{\sgc}{{\mathfrak c}}
\newcommand{\sgd}{{\mathfrak d}}

\newcommand{\sC}{{\cal C}}

\newcommand{\sF}{{\cal F}}

\newcommand{\sK}{{\kappa}}

\newcommand{\sS}{{\cal S}}

\newcommand{\sZ}{{\cal Z}}

\newcommand{\gaa}{{\gamma}}

\newcommand{\cw}{w}
\newcommand{\tc}{{\tilde{c}}}

\newcommand{\zk}{Z_K}
\newcommand{\zq}{Z_\QQ}
\newcommand{\hzk}{\hat{\zeta}_K}

\newcommand{\ttc}{{c^{'}}}
\newcommand{\txi}{{\xi}}

\newcommand{\xiq}{\xi_\QQ}

\newcommand{\eqn}[1]{(\ref{#1})}

\makeatletter
\def\eqalignno#1{\displ@y \ta {\bf s} kip\@centering
  \halign to\displaywidth{\hfil$\@lign\displaystyle{##}$\ta {\bf s} kip\z@skip
    & $\@lign\displaystyle{{}##}$\hfil\ta {\bf s} kip\@centering
    & \llap{$\@lign##$}\ta {\bf s} kip\z@skip\crcr
    #1\crcr}}
\makeatother

\makeatletter
\def\@sect#1#2#3#4#5#6[#7]#8{\ifnum #2>\c@secnumdepth
     \def\@svsec{}\else 
     \refstepcounter{#1}\edef\@svsec{\csname the#1\endcsname.\hskip .75em }\fi
     \@tempskipa #5\relax
      \ifdim \@tempskipa>\z@ 
        \begingroup #6\relax
          \@hangfrom{\hskip #3\relax\@svsec}{\interlinepenalty \@M #8\par}%
        \endgroup
       \csname #1mark\endcsname{#7}\addcontentsline
         {toc}{#1}{\ifnum #2>\c@secnumdepth \else
                      \protect\numberline{\csname the#1\endcsname}\fi
                    #7}\else
        \def\@svsechd{#6\hskip #3\@svsec #8\csname #1mark\endcsname
                      {#7}\addcontentsline
                           {toc}{#1}{\ifnum #2>\c@secnumdepth \else
                             \protect\numberline{\csname the#1\endcsname}\fi
                       #7}}\fi
     \@xsect{#5}}
\def\@begintheorem#1#2{\it \trivlist \item[\hskip \labelsep{\bf #1\ #2.}]}
\makeatother

\numberwithin{equation}{section}

\begin{document}
\begin{center}
{\Large {\bf On a Two-Variable Zeta Function for Number Fields }} \\
\vspace{\baselineskip}
{\em Jeffrey C. Lagarias}~\footnote{Work done in part during a visit
to the Institute for Advanced Study.}  \\
{\em Eric Rains} \\
\vspace*{.2\baselineskip}
AT\&T Labs - Research \\
Florham Park, NJ 07932 \\

\vspace*{1\baselineskip}
(July 6, 2002 version) \\
\vspace{2\baselineskip}
{\em Abstract}
\end{center}
This paper studies  a two-variable zeta
function  $\zk (\cw, s)$ attached to an algebraic number field
$K$, introduced by  van der Geer and Schoof~\cite{vdGS99}, which is 
based on an
analogue of the Riemann-Roch theorem for number fields using
Arakelov divisors. 
When $\cw= 1$ this function becomes the completed Dedekind zeta 
function $\hzk(s)$
of the field $K$.  The function
is an meromorphic  function of
two complex variables with polar divisor $s(\cw - s)$, 
and it  satisfies the functional equation
$\zk(\cw, s) = \zk(\cw, \cw - s)$.
 We consider the special case $K = \QQ$,
where for $\cw = 1$ this function 
is $\hat{\zeta}(s)= \pi^{-\frac{s}{2}} \Gamma(\frac{s}{2}) \zeta(s)$. 
The function $\xi_{\QQ}(\cw, s) := \frac{s(s- \cw)}{2\cw} \zq(\cw, s)$ is 
shown to be an entire function on
$\CC^2$, to satisfy the functional equation
$\xi_{\QQ}(\cw, s) = \xi_{\QQ}(\cw, \cw - s),$ and
to have 
$\xi_{\QQ}(0, s) =-\frac{s^2}{8}(1 - 2^{1 + \frac{s}{2}})
(1 - 2^{1 - \frac{s}{2}}) \hat{\zeta}(\frac{s}{2}) 
\hat{\zeta}(\frac{-s}{2}).$
We study the location of the zeros of $\zq(\cw, s)$ 
for various real values
of $\cw = u$. For fixed $u \geq 0$
the zeros are confined to a vertical strip of width
at most $u + 16 $ and the number of zeros $N_u(T)$ to height $T$ has
similar asymptotics to the Riemann zeta function.
For fixed $u < 0$ 
these functions are strictly positive on the ``critical line''
$\Re(s) = \frac{u}{2}$. 
This phenomenon is associated to a positive convolution semigroup with
parameter $u \in \RR_{> 0}$, which is a semigroup
of infinitely divisible probability distributions.
having densities $P_u(x)dx$ for real $x$, where 
$P_u(x) = \frac{1}{2\pi}\theta(1)^u Z_{\QQ}(-u, -\frac{u}{2} + ix),$
and $\theta(1) = \pi^{1/4}/ \Gamma(3/4).$ \\



\setlength{\baselineskip}{1.0\baselineskip}

%
%
%

\section{Introduction}
\hsp
Recently van~der Geer and Schoof \cite[Prop. 1]{vdGS99} 
formulated an
``exact'' analogue of the Riemann-Roch theorem valid  for an algebraic
number field $K$, based on Arakelov divisors.
They used  this result to formally express the completed zeta function
$\hzk (s)$ of $K$ as an integral over the
 Arakelov divisor class group $Pic(K)$ of $K$. They introduced
a two-variable zeta function attached to  a number field $K$, 
also given as an integral
over the Arakelov class group, which
we call either the {\em Arakelov zeta function}
or the {\em two-variable zeta function}.
This zeta function was
modelled after a two-variable zeta function attached to 
a function field over a finite field, introduced
in 1996 by Pellikaan~\cite{Pe96}. For convenience
we review the Arakelov divisor interpretation
of the two-variable zeta function and the Riemann-Roch theorem for
number fields in an appendix.

In this  paper we study in detail the 
two-variable zeta function 
attached to the rational field $K = \QQ$.
Then in the final section we consider   
two-variable zeta functions for general
algebraic number fields $K$.
The results are derived starting from an integral representation of
this function, and if one takes it as given, then
the paper is independent of the Arakelov divisor interpretation.
The Arakelov class group of $\QQ$ can be identified
with the positive real line (with multiplication as the group
operation) and 
van der Geer and Schoof's integral
becomes, formally,  
\begin{equation}\label{S103}
\zq (s) \cong \int_0^\infty \theta (t^2)^s \theta(\frac{1}{t^2})^{1 - s} 
\frac{dt}{t},
\end{equation}
in which
\begin{equation}\label{S104}
\theta (t) := \sum_{n \in \ZZ} e^{- \pi n^2 t}
\end{equation}
is the theta function  $\theta(t) = \vartheta_3 (0, e^{\pi t})$, where
$$\vartheta_3 (z, q) := \sum_{n \in \ZZ} e^{2 \pi i n z} q^{n^2}
= 1 + 2\sum_{n = 1}^\infty q^{n^2} \cos 2nz, $$
is a Jacobi theta function. 
The $\cong$ used in \eqn{S103} reflects the fact that the integral 
on its right side converges {\em nowhere};
a regularization is needed to assign it a meaning. 
Such a regularization can be obtained using the 
Arakelov two-variable zeta function $\zq (\cw,s)$ attached to $\QQ$,
which we define to be
\begin{equation}~\label{S105}
\zq (\cw,s)  := \int_0^\infty \theta (t^2 )^s \theta (\frac{1}{t^2})^{\cw-s}
\frac{dt}{t} ~.
\end{equation}
Our definition here differs from the one
 in van der Geer and Schoof \cite{vdGS99}
by a linear change of variable, setting their second variable $t = \cw - s.$
The integral on the right side of \eqn{S105}
has a region of absolute
convergence in  $\CC^2$, which is the open cone
\begin{equation}~\label{105b}
\sC:= \{(\cw, s): \Re(\cw) < \Re(s) < 0 \}.
\end{equation}
The function $\zq (\cw,s)$ meromorphically continues from
the cone $\sC$ to all of
$\CC^2$, with polar divisor consisting of the (complex) hyperplanes
$\{s=\cw\} \cup \{s=0\}$, a set of real-codimension two, see \S2.
On restricting $\zq (\cw,s)$  to the line $\cw=1$,
the resulting function is the {\em completed Riemann zeta function} 
$\hat{\zeta}(s)$, which is
$$
\hat{\zeta} (s) := \pi^{- s/2} \Gamma (s/2) \zeta (s).
$$
Thus the two-variable zeta function $\zq (\cw,s)$ 
defined via \eqn{S105}  provides a 
method to regularize the integral \eqn{S103}, 
and the same can be done for arbitrary number fields.

We are motivated by several questions about this function.

(1) What are the properties of the function as a meromorphic
function of two complex variables? In particular, determine 
information about its zero divisor.

(2) What is the meaning of the additional variable $\cw$ and what
arithmetic information does it encode?

(3) What properties of this two-variable zeta
function reflect Arakelov geometry?

(4) Is there any connection between zeta functions encoding
information based on Arakelov geometry and zeta functions coming
from automorphic representations and the
Langlands program?

This paper mainly addresses question (1),
obtaining information on the zero set of the Arakelov zeta function.
Concerning question (2), we observe in \S2 that
the function  $\zq(\cw, s)$ is representable by the integral
\begin{equation}~\label{107a}
\zq (\cw,s) = \frac{1}{2} 
\int_{0}^\infty \theta(t)^{\cw} t^{s/2} \frac{dt}{t},
\end{equation}
which expresses it as a Mellin transform of 
$2\theta(t)^{\cw}$. The function $\theta(t)^{\cw}$ is a modular
form of weight $\frac{\cw}{2}$ (with multiplier system)
on a congruence subgroup of
the modular group, and the complex
variable $\cw/2$ parametrizes the {\em weight} of this modular form.
The arithmetic information it encodes includes the
invariants $\eta(\QQ)$ and $g(\QQ)$ introduced in 
van der Geer and Schoof~\cite{vdGS99}, defined in the appendix.
Concerning question (3), we observe that there is an
extra structure  associated to $\zq(\cw, s)$, which
is a holomorphic convolution semigroup of
complex-valued measures on lines 
$L(t) = \{(u, u+it): -\infty < t < \infty \}$  in the 
real-codimension one cone
$$
\sC^{-} := \{(\cw, s): \cw= u \in \RR ~~\mbox{and}~~ u < \Re(s) < 0 \}\,
$$
see \S7.2.
This cone is contained in the region of absolute convergence $\sC$
of the integral representation \eqn{S105}. Of
particular interest for $\zq(\cw, s)$ is the real-codimension two subcone 
$$
\sC_{crit} = \{ (\cw, s) : \cw = u \in \RR~~\mbox{and}~ \Re(s) = \frac{u}{2} < 0\},
$$
which generalizes the ``critical line'' of the zeta function,
and on which the measures are real-valued.
Perhaps this semigroup structure is  associated in some way
with Arakelov geometry, since various constants associated with the
semigroup on the subcone $\sC_{crit}$ have  arithmetic interpretations in the
framework of van der Geer and Schoof, see  \S7. 
Concerning question (4), the subject of Arakelov geometry was
developed in part to answer Diophantine questions and 
has a completely different origin from automorphic representations.
Any connection between these two subjects
could potentially be of great interest.
However we do not find any obvious connection, and note only that
the $\cw$ variable
interpolates between modular forms of different weights,
and when $\cw$ is a positive even integer 
these are holomorphic modular
forms of the type appearing in automorphic representations.
In general these forms are not eigenforms for 
Hecke operators, and in \S3.4 we show these forms have 
associated Euler products
exactly when $\cw= 0, 1, 2, 4$ and $8$.

Besides giving information on questions (1)-(4) above,
the analysis of this paper may be useful for other purposes. 
This function provides an interesting example of an
entire function in two complex variables of finite order,
see Ronkin~\cite{Ro89} and Stoll~\cite{St74}.
The information about the zero
locus of $\zq(\cw, s)$ that we obtain mainly concerns the region where
where the variable $\cw$ is real; these  are Mellin transforms
of modular forms of real weight, which have been
extensively studied.  The movements of zeros in the $s$-plane
as the (real) parameter $\cw$ is varied  may be compared 
with movement of zeros under milder deformations such as
those in linear combinations of L-functions, see
Bombieri and Hejhal~\cite{BH95} and Hejhal~\cite{He00}.

The function $\zq (\cw,s)$ shares many properties of 
the Riemann zeta function. It satisfies the functional equation
\begin{equation}\label{S106}
\zq (\cw,s) = \zq (\cw,\cw-s) ~.
\end{equation}
When $\cw= u$ is real then $\zq (u,s)$ retains several familiar symmetries of 
the Riemann zeta function: it is real on the real axis 
$\Im (s) =0$, and it is real on the
``critical line'' ${\rm Re} (s) = \frac{u}{2}$, 
which is the line of symmetry of the functional equation.
Thus for fixed real $u$, the zeros of 
$\zq (u,s)$ which do not lie on the critical line or 
the real axis must occur in sets of four:
$s$, $u-s$, $\bar{s}$, $u- \bar{s}$.  
This extra symmetry can be used to extract information about 
the components of the zero locus, see Lemma~\ref{Nle71}.
On the other hand, for most real $u$ the function
$\zq (u,s)$ fails to satisfy a Riemann hypothesis,
as we describe below. 
The  Riemann zeta function appears when $\cw = 1$ and it is
interesting to note that it also appears 
in terms of data at $\cw =0$, given
in the following result, which is proved in \S5.

\begin{theorem}~\label{th10}
The function
$$\txi_{\QQ}(\cw, s) = \frac{s(s-\cw)}{2\cw}Z_{\QQ}(\cw, s)$$ 
is an
entire function in two complex variables. At $\cw=0$ it is  
\begin{equation}~\label{107}
\txi_{\QQ}(0, s)= -\frac{s^2}{8}(1 - 2^{1 + \frac{s}{2}})
(1 - 2^{1 - \frac{s}{2}}) \hat{\zeta}(\frac{s}{2}) 
\hat{\zeta}(\frac{-s}{2}),
\end{equation}
where $\hat{\zeta}(s)$ is the completed Riemann zeta function.
In particular, for all real $t$,
\begin{equation}\label{108a}
\txi_{\QQ}(0, it) = \frac{t^2}{8} 
|1 - 2^{ 1 + \frac{it}{2}}|^2 
|\hat{\zeta}(\frac{it}{2} )|^2
\end{equation}
is strictly positive, with  $\txi_{\QQ}(0, 0) = 1/2.$
\end{theorem}

The Jacobi triple product formula plays an essential role in
our derivation
of the formula \eqn{107}. Note that the function
$\txi_{\QQ}(1,s)$ coincides with the Riemann $\xi$-function,
and  the functional equation 
$\txi_{\QQ}(\cw, s)=  \txi_{\QQ}(\cw,\cw -s)$ is inherited
from $Z_{\QQ}(\cw, s).$

The most striking result of this paper appears in \S7,
and concerns for negative real $\cw= -u~ (u > 0)$ the behavior
of the function $\zq(\cw, s)$ on the
``critical line'' $\Re(s) = - \frac{u}{2}$. This result
makes a connection with probability theory, involving  infinitely
divisible distributions.

\begin{theorem}~\label{th11}
For negative real $\cw$, with $\cw =-u$  $(u > 0)$ the function 
$\zq\left( -u, -\frac{u}{2} + it \right)$ is given
as the Fourier transform 
\begin{equation}~\label{105a}
 \zq\left( -u, -\frac{u}{2} + it \right) = 
\theta(1)^{-|u|} \int_{-\infty}^\infty
f(r)^{u} e^{-irt} dr
\end{equation}
in which
$$f(r) = \theta (1)\frac{e^{r/2}}{\theta (e^{-2r} )} 
=\frac{\theta (1)}{\sqrt{\theta (e^{-2r}) \theta (e^{2r} )}} 
= \theta (1)\frac{e^{-r/2}}{\theta (e^{2r} )}~$$ 
and
$$
\theta(e^{-2r}) = 1 + 2\sum_{n=1}^\infty e^{-\pi n^2 e^{-2r}}.
$$
The function  $f(r)$
is the characteristic function of an infinitely
divisible probability measure with finite second moment,
whose associated Khintchine canonical measure
$\frac{1}{1+x^2}M\{dx\}$ has  $M \{dx\}= M(x) dx$ with
\begin{equation}\label{105}
M(x) = \frac{1}{8\pi}x^2
\left| \hat{\zeta} \left( \frac{ix}{2} \right) \right|^2
| 1 - 2^{1 + \frac{ix}{2}}|^2 ~ = \frac{1}{\pi} \xi_{\QQ}(0, ix),
\end{equation}
in which $\hat{\zeta} (s)$
is the completed Riemann zeta function. In particular,  
\begin{equation}\label{105aa}
\zq\left( -u, -\frac{u}{2} + it \right) > 0 
\qquad\mbox{for}\qquad -\infty < t < \infty.
\end{equation}
\end{theorem}

Many  interesting connections
between zeta and theta functions and probability theory are known;
see Biane, Pitman and Yor~\cite{BPY99} for a comprehensive survey.
Theorem~\ref{th11} appears structurally different from
any of the known results. 
The  positivity property
 \eqn{105aa} can be called an ``anti-Riemann hypothesis'',
because it shows there are no zeros on the ``critical line''
${\rm Re} (s) = -\frac{u}{2}$ for fixed real $u > 0$.

There are a number of different canonical forms used to specify
infinitely divisible distributions.
Feller~\cite[pp. 563]{Fe71} uses
the canonical measure $M\{dx\}$,
which we term the 
{\em Feller canonical measure},
while the  {\em Khintchine canonical measure} 
$K\{dx\}   =\frac{1}{1+x^2}M\{dx\}$ is often used, 
see \cite[pp.564-5]{Fe71}. 
An infinitely divisible distribution is a member of a positive
convolution semigroup of measures, and the Feller canonical
measure is related to the infinitesimal generator of the
semigroup. The measure  $M(x)dx$ above
involves the values of the Riemann zeta function on the boundary of
its critical strip, noting that the functional equation gives
$$ \hat{\zeta}(\frac{ix}{2}) = \hat{\zeta}(1 - \frac{ix}{2}).$$
 Theorem~\ref{th11} follows from two results proved in \S7,
Theorem~\ref{nth61} and  Theorem~\ref{nth62}.
We also note that the value
$\theta(1)= \frac{\pi^{\frac{1}{4}}}{\Gamma \left( \frac{3}{4} \right)}
\approx 1.08643$ appearing in Theorem~\ref{th11}
equals $e^g$ where $g$ is the ``genus of $\QQ$''
as defined by van der Geer and Schoof~\cite{vdGS99}, see the appendix.

The positive holomorphic convolution semigroup structure
associated to this two-variable zeta function merits further study.
It seems an interesting  question to determine the generality of this
positivity property.
All algebraic number fields $K$ have an associated holomorphic
convolution semigroup of complex-valued
measures, which are real-valued measures on the ``critical line''.
However the positivity fails to hold in general, and
perhaps is true only for a few specific
number fields, see \S9.

We comment on related work. There is precedent for studying 
two-variable functions given
by integrals of the form (\ref{107a}) with $\theta(t)$ replaced with
some other modular form. Conrey and Ghosh \cite[Sect. 5]{CG94} 
considered a Fourier integral associated to powers of the modular form  
$\Delta(\tau)$ of weight 12, which is a cusp form.
The integral they consider
can be transformed to a constant multiple of the integral  
\begin{equation}~\label{111a}
\tilde{Z}(\cw, s)= \int_{0}^\infty \Delta(it)^{\cw/24} t^{s/2} \frac{dt}{t},
\end{equation}
where they take $\cw= k >0$, and $s=it$. 
They note that  associated Dirichlet series has an Euler product for
$\cw = 1, 2, 3, 4, 6, 8, 12,$ and $24$. 
Bruggeman \cite{Br94} studied properties of 
families of automorphic forms 
of variable weight; he considers powers of the Dedekind  eta
function family $\eta(t)^{\cw}$ 
in \cite[1.5.5]{Br94}.
This family appears in 
\eqn{111a} since $\Delta(\tau) = \eta(\tau)^{24}$,
see \cite[p. 11]{Br94}.

We now summarize the contents of the paper.
In \S2 we give the analytic continuation and functional equation
for $\zq(\cw, s),$ essentially following Riemann's second proof
of the functional equation for $\zeta(s)$. We derive integral
formulas for $Z_{\QQ}(\cw, s),$ which converge on 
$\CC \times \CC$ off certain hyperplanes.

In \S3,
as a preliminary to later results,
we study the Fourier coefficients $c_m(\cw)$ of the modular form
$$\theta(it)^{\cw} = 1 + \sum_{m=1}^\infty c_m(\cw) e^{-\pi i m t}.$$ 
We show that  $(-1)^m m! c_m(-\cw)$ is a polynomial of degree $m$
with nonnegative integer coefficients. For $\cw= u$ on the positive real
axis we obtain the estimate  
$|c_m(u)| \le 6u m^{\frac{u}{2}+1},$
 whose merit is that it is uniform in $u$.
For general  $|\cw| = R$ there is an  upper bound
$$|c_m(\cw)| \leq C_0 R^{\frac{R}{2}}e^{\pi \sqrt{Rm}},$$
which follows from classical estimates. 
We show that
the Dirichlet series $\tilde{D}_{\cw}(s)= \frac{1}{2w} D_{\cw}(s)$
for $\cw \in \CC$ has an Euler product if and only if
$\cw= 0, 1, 2, 4$ and $8$.

In \S4 we study growth properties of the entire function
$\xi_{\QQ}(\cw, s) := \frac{s(s-\cw)}{2\cw}\zq(\cw, s).$
We first show that $\xi_{\QQ}(\cw, s)$  is an 
entire function of order one
and infinite type in two complex variables, in the sense that
it satisfies the growth bound:
There is a constant $C_1$ such that for any $(\cw, s) \in \CC^2$,
if  $R = |s| + |\cw| + 1$, then 
$$|\xi_{\QQ}(\cw, s)| \leq e^{C_1 R \log R}.$$
 Thus any
linear slice function $f(s) = \xi_{\QQ}(as + b, cs + d)$,
has at most 
$O(R \log R)$ zeros in the disk of radius $R$ as $R \to \infty$, provided it 
is not identically zero.
We then show that for fixed $\cw \in \CC$ and fixed 
$\sigma \in \RR$, the function 
$$f_{\cw, \sigma}(t) := \xi_{\QQ}(\cw, \sigma + it)~~~-\infty < t < \infty,$$ 
has rapid decrease,
is in the Schwartz class $\sS(\RR),$ and 
is uniformly bounded in vertical strips
$\sigma_1 \leq \sigma \leq \sigma_2$, for finite $\sigma_1, \sigma_2$.

In \S5 we treat the case $\cw=0$ and prove Theorem~\ref{th10}.

In \S6 we treat the case when $\cw=u >0 $ 
is a fixed positive real number, and study the 
zeros of $\xi_{\QQ}(u, s)$. 
We show that these zeros are confined
to the vertical strip
$|\Re(s) - \frac{u}{2}| < \frac{u}{2} + 8.$
Then we show that the number $N_u(T)$ of zeros $\rho$ having
$|\Im(\rho)| \leq T$ has similar asymptotics to that of the
Riemann zeta function, namely
$$\frac{1}{2} N_u (T) = \frac{T}{2 \pi} \log \frac{T}{2 \pi} 
- \frac{T}{2 \pi} + \frac{7}{8} + S_u (T),$$
with
$$|S_u (T) | \le C_0 (u + 1) \log (T+ u + 2).$$
and the constant $C_0$ is absolute.
The zeros of $\xi_{\QQ}(u, s)$ appear
to lie on the ``critical line'' $\Re(s) = \frac{u}{2}$ only for special
values $u =1$ and $u=2$;
we observe that only an infinitesimal fraction
of zeros are on this line for $u= 4$ and $u = 8$.

In \S7 we consider $\xi_{\QQ}(\cw, s)$ where $\cw=-u$ is a fixed
negative real number ($-u < 0$). 
We prove Theorem ~\ref{th11}, that the function
$\xi_{\QQ}(-u, s)$, which is necessarily real on the critical
line  $\Re(s) = -\frac{u}{2}$, is always positive there.
The proof of this result makes essential use of the
Jacobi product formula, which is applicable because the constant
term in the theta function is present in the integral representation
\eqn{107a}. The associated structure behind Theorem~\ref{th11}
is a holomorphic convolution
semigroup $\rho_{u,v}(x) dx $ of complex-valued measures on
the real line, defined for $(u, v)$ real in the cone
$u > 0$ and $|v| < u$, and these measures are
positive real on the line $v = 0$. We derive formulae for the
cumulants and moments of these measures.
We also list a number of open questions concerning the
location of zeros for negative real $u$.
For example, for real $\cw= u$, 
are the asympotics of the number
of zeros $\rho$ with $|\Im(\rho)| < T$ as $T \to \infty$
the same for negative real $u$ as they are for positive real $u$?

In \S8 we consider general complex $\cw$, and the 
zero locus $\sZ_{\QQ}$  of $\xi_{\QQ}(\cw, s)$. 
The set $\sZ_\QQ$ viewed 
geometrically~\footnote{Considered algebraically there is the 
additional problem of determining the multiplicity of each 
irreducible component.}
is a one-dimensional complex manifold, having more than one 
irreducible analytic component (possibly infinitely many components), 
each one of which is a Riemann surface embedded in $\CC^2$.
We show that the zeta zeros
$\rho_7$ and $\rho_8$ are on the same irreducible component, and raise
the question whether the zeta zeros (for $\cw=1$, $s$ varying) are 
all on a single irreducible component of the zero locus.

In \S9  we briefly consider Arakelov
zeta functions attached to general algebraic  number fields $K$.
All results of this paper extend to the Arakelov
zeta function attached to the Gaussian field $K = \QQ(i),$
and many of the results
extend to general $K$, with similar proofs. However 
the positivity property of
Theorem~\ref{th11} for $K=\QQ$, our proof used a product formula 
for the modular form and does not extend
to general number fields $K$. Numerical experiments 
show that the positivity property does not hold for 
several imaginary quadratic fields, with discriminants $-8, -11$ and
$-19$. Our computations
 allows the possibility that it might hold for some fields
whose modular forms do not have a product formula, including
the imaginary quadratic fields with discriminants $-3$ and $-7$.

In the appendix we review the Arakelov divisor framework of van der Geer
and Schoof~\cite{vdGS99}, and derive formulas for the two-variable
zeta function for $K= \QQ$ and $\QQ(i)$.

\noindent\paragraph{Acknowlegments.} The authors thank
E. Bombieri, J. B. Conrey, C. Deninger, J. Pitman, J. A. Reeds, 
M. Yor and the reviewer for helpful comments and references. \\

\noindent \paragraph{Notation.}
The variables $w, s, z$ denote  complex variables with
$w= u + iv$, $s= \sigma + it$, $z= x+iy$, and $u, v, \sigma, t, x, y$
always denote real variables.
We use two  versions of the Fourier transform, 
differing in their scaling, because the usual conventions for
the Fourier transform differ in probability
theory and number theory.
 The {\em Fourier transform} $\sF $ is given by
$$ \sF f (x) : = \int_{-\infty}^\infty f(t) e^{-2\pi i xt} dt, $$
with inverse
$$ \sF^{-1} f(t) :=  \int_{-\infty}^\infty f(x) e^{2\pi i xt} dx. $$
In probability theory the {\em characteristic function} 
$\varphi = \varphi_M$ of a
Borel measure $M\{dx\}$ of unit mass on the line is 
$$\varphi_M(t) := \int_{-\infty}^\infty  e^{ixt}M\{dx\}, $$
In the case where $M\{dx\}= f(x) dx$ we write $\varphi_M(t) = \check{f}(t)$
as an inverse Fourier transform, and 
the corresponding Fourier transform is
\begin{equation}~\label{111}
\hat{f}(x) = \frac{1}{2\pi}  \int_{-\infty}^\infty f(t) e^{-i xt} dt.
\end{equation}
%
%

\section{Analytic Continuation and Functional Equation}
\hsp
We now  obtain the meromorphic continuation of $\zq(\cw, s)$, 
which determines its polar divisor and part of its zero divisor.
Using the theta function transformation formula
\begin{equation}\label{S202}
\theta (t^2) = \frac{1}{t} \theta \left( \frac{1}{t^2} \right) 
\quad\mbox{for}\quad t > 0
\end{equation}
we  can rewrite
\begin{equation}\label{S201}
\zq (\cw,s) = \int_0^\infty \theta \left( t^2 \right)^s \theta 
\left( \frac{1}{t^2} \right)^{\cw-s} \frac{dt}{t}
\end{equation}
in the form
\begin{equation}\label{S203}
\zq (\cw,s) = \int_0^\infty \theta \left( \frac{1}{t^2} \right)^{\cw}
t^{-s} \frac{dt}{t}.
\end{equation}
Then, after a change of variable $t \to 1/t$, followed by $t^2 \to u$, 
one obtains
\begin{equation}\label{S204}
\zq (\cw,s) = \int_0^\infty \theta (t^2)^{\cw} t^s \frac{dt}{t}
= \frac{1}{2}\int_0^\infty \theta (u)^{\cw} u^{s/2} \frac{du}{u}.
\end{equation}
Note that
$\theta (t^2) -1 \to 0$ rapidly as $t \to \infty$, hence
$\theta (t^2) - \frac{1}{t} \to 0$ rapidly as $t \to 0^+$.
This implies that \eqn{S201} converges absolutely on the open domain
$\sC = \{ (\cw,s) : \Re(\cw) < \Re (s) < 0 \}$ in $\CC^2$.
The convergence is uniform on compact subsets of this domain, which defines
$\zq (\cw,s)$ as an analytic function there.

\begin{theorem}\label{th21}
The function $\xi_\QQ (\cw,s) = \frac{s(s - \cw)}{2\cw} \zq (\cw,s)$ 
analytically 
continues to an entire function on $\CC^2$, and 
satisfies the functional equation
\begin{equation}\label{S205}
\xi_\QQ (\cw,s) = \xi_\QQ (\cw, \cw-s) ~.
\end{equation}
\end{theorem}

\noindent\paragraph{Remark.}
For $\cw=1$ we have
$\zq (1,s) = \pi^{-s/2} \Gamma ( s/2) \zeta (s)$ and
$\xi_\QQ (1,s) = \xi (s)$ where
$\xi(s): = \frac{1}{2} s(s-1) \pi^{-s/2} \Gamma (s/2) \zeta (s)$ 
is Riemann's $\xi$-function, and we recover the functional
equation for $\zeta(s)$. We give a Fourier-Laplace transform
integral representation for
$\xi_\QQ (\cw,s)$ in Theorem~\ref{th44}.

\paragraph{Proof.}
We split the integral \eqn{S201} into two pieces $\int_0^1$ 
and $\int_1^\infty$ and consider them separately.
Using the transformation law yields
\begin{eqnarray}\label{S206}
 \int_0^1 \theta (t^2)^s \theta \left ( \frac{1}{t^2} \right)^{\cw-s}
\frac{dt}{t} & = & \int_0^1 \theta \left( \frac{1}{t^2} \right)^{\cw} t^{-s} 
\frac{dt}{t} \nonumber \\
&= & \int_0^1 \left( \theta \left( \frac{1}{t^2} \right)^{\cw} - 1 \right)
t^{-s} \frac{dt}{t} + \int_0^1 t^{-s} \frac{dt}{t} \,.
\end{eqnarray}
Both sides are defined and converge ${\rm Re} (s) < 0$ and 
${\rm Re} (\cw) < 0$.
On the right side the first integral converges for all 
$(\cw,s) \in \CC \times \CC$, because for $|\cw| \le R$ and $0 \le t \le 1$,
\begin{equation}\label{S207}
\theta \left( \frac{1}{t^2} \right)^{\cw} = 1 + \cw e^{- \frac{\pi}{t^2}} +
O \left( e^{- \frac{(2 - \epsilon ) \pi}{t^2}} \right)
\end{equation}
as $t \to 0^+$, where the constant in the O-symbol depends only on $R$.
This uniformity of convergence shows that this integral 
is an entire function on $\CC^2$.
The second integral in \eqn{S206} converges, for 
${\rm Re} (s) < 0$, to the function $- \frac{1}{s}$.

Similarly,
\begin{eqnarray}\label{S208}
\int_1^\infty \theta (t^2)^s \theta \left( \frac{1}{t^2} \right)^{\cw -s }
\frac{dt}{t} & = &
\int_0^1 \theta \left( \frac{1}{t^2} \right)^s \theta (t^2)^{\cw - s}
\frac{dt}{t} \nonumber \\
& = & \int_0^1 \theta \left( \frac{1}{t^2} \right)^{\cw} t^{s-\cw} 
\frac{dt}{t} \nonumber \\
& = & \int_0^1 \left( \theta \left( \frac{1}{t^2} \right)^{\cw} -1 \right)
t^{s-\cw} \frac{dt}{t} + \int_0^1 t^{s-\cw} \frac{dt}{t} \,,
\end{eqnarray}
with both sides convergent for ${\rm Re} (s) < 0$, 
${\rm Re} (\cw) < 0$ and ${\rm Re} (s - \cw) > 0$.
This region overlaps the region of convergence of \eqn{S204} 
in an open domain in $\CC^2$.
The first integral on the right side of \eqn{S207} 
defines an entire function on $\CC^2$, while 
the second integral in \eqn{S208} is $\frac{1}{s- \cw}$
for ${\rm Re} (\cw -s) > 0$.
We obtain
\begin{equation}\label{S209}
\zq (\cw,s) = - \frac{1}{s} +
\frac{1}{s-\cw} + \int_0^1
\left( \theta \left( \frac{1}{t^2} \right)^{\cw} -1 \right)
(t^{-s} + t^{-(\cw-s)} )
\frac{dt}{t}
\end{equation}
which is valid for $(\cw,s) \in \CC^2$, for $s \neq 0, \cw$.
Since the right side of this equation is invariant under 
$s \to \cw -s$, we obtain the functional  equation
\begin{equation}\label{S210}
\zq (\cw,s) = \zq(\cw,\cw -s) ~.
\end{equation}
Now \eqn{S209} implies that $s(\cw -s) Z_{\QQ}(\cw, s)$ is
an entire function on $\CC^2$,
and on  setting $\cw = 0$ in  \eqn{S209} we see that
$Z(0, s)$ is identically zero. Thus 
$\txi_{\QQ}(\cw, s) = \frac{s(s - \cw)}{2\cw}Z_{\QQ}(\cw, s)$ is
also an entire function on $\CC^2$, and satisfies the 
same functional equation.~~~$\bsq$

Viewed as a modular form, $\theta (t)$ in  \eqn{S201} is not a cusp form, 
due to its nonzero constant term.
One consequence is that the Mellin transform
$\int_0^\infty \theta (t) t^s \frac{dt}{t}$ fails to converge anywhere.
Riemann's second proof of the functional equation (for $\cw=1$)
circumvents 
this problem by removing the constant term, using 
$2 \psi (t) = \theta (t) -1$ in the integrand, and in this case  
the Mellin transform integral converges for $\Re(s) > 1$.
In Theorem \ref{th21}, the constant term 
``evaporates'' because, formally,
\begin{equation}\label{S211}
\int_0^\infty t^s \frac{dt}{t} \equiv 0 ~.
\end{equation}
More precisely
\begin{eqnarray*}
\int_0^1 t^s \frac{dt}{t} & = & \frac{1}{s} \quad\mbox{for}\quad
{\rm Re} (s) > 0 \,, \\
\int_1^\infty t^s \frac{dt}{t} & = & - \frac{1}{s}
\quad\mbox{for}\quad {\rm Re} (s) < 0 \,.
\end{eqnarray*}
One convention for ``regularization''
of the integral is to analytically continue these two pieces separately 
and then add them, which results in \eqn{S211}.
Theorem \ref{th21} justifies  this convention
by introducing the extra variable $\cw$, 
finding a common domain $\sC$ in the $(\cw,s)$-plane 
where the integral converges, and then analytically continuing 
in both variables to the line $\cw=1$.

We next give  modified integral formulas for $\zq(\cw, s)$
valid on most of  $\CC^2$.
We define {\em Heaviside's function} $H(s)$ for complex $s$
with $\Re(s) \neq 0$ to be
\[
H(s) = H(\Re(s)) = \begin{cases} 1 & \Re(s)>0\\
                                \frac{1}{2} & \Re(s)=0\\      
                                 0 & \Re(s)<0\end{cases}.
\]
\begin{theorem}~\label{th22a}
If $\Re(s)\notin \{\Re(\cw),0\}$, then
\begin{align}
Z_\QQ(\cw ,s) 
&= \int_0^\infty (\theta(t^2)^{\cw}-H(s)-H(\cw-s)t^{-\cw}) t^s \frac{dt}{t}, 
\nonumber \\
&= \int_{-\infty}^\infty (\theta(e^{2x})^{\cw}-H(s)-H(\cw-s)e^{-\cw x}) 
e^{sx} dx,
\end{align}
where both integrals converge absolutely.
\end{theorem}

\begin{figure}[hbp]~
\begin{center}
\input{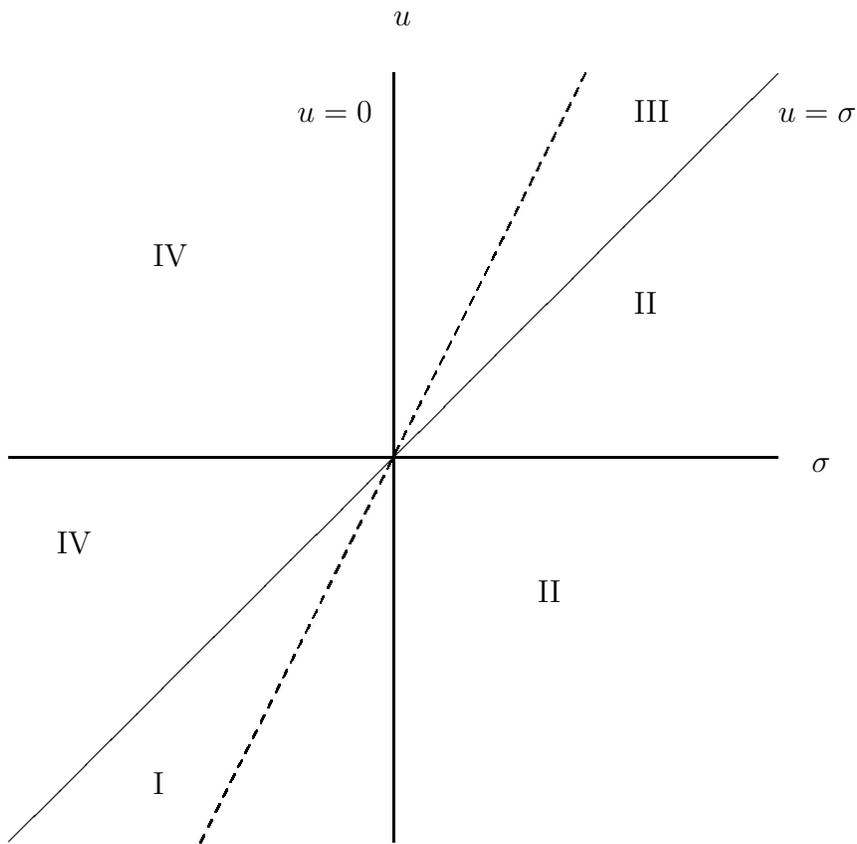}
\end{center}
\caption{Convergence regions: $u = \Re(\cw), \sigma = \Re(s).$}
\label{fig21}
\end{figure}

\noindent\paragraph{Remark.} This result expresses  $Z_{\QQ}(\cw,s)$
by a convergent  integral formula with integrand
\begin{align}
(I)~~~~~~~~~~~~~~~~~~ \theta(t^2)^{\cw}  
& \qquad\mbox{if}\qquad  \Re(\cw) < \Re(s) < 0 \, \\
(II)~~~~~~~~~~~~ \theta(t^2)^{\cw} - 1 
& \qquad\mbox{if}\qquad \Re(\cw - s) < 0 < \Re(s) \, \\
(III)~~~ \theta(t^2)^{\cw} - 1 - t^{-\cw} 
& \qquad\mbox{if}\qquad 0 < \Re(s) < \Re(\cw) \, \\
(IV)~~~~~~~~~ \theta(t^2)^{\cw} - t^{-\cw}
& \qquad\mbox{if}\qquad  \Re(s) < 0  < \Re(\cw - s)  \, 
\end{align}
These regions are pictured in Figure~\ref{fig21}; the dotted line
is the ``critical line'' $\sigma = \frac{\Re(\cw)}{2}.$

\paragraph{Proof.}
The second integral (Laplace transform) follows by the
change of variable $t = e^x$, so it suffices to consider the
first integral.
We recall
\begin{align}
Z_\QQ(\cw ,s)
&=
-\frac{1}{s}+\frac{1}{s-\cw} + 
\int_0^1
(\theta(1/t^2)^{\cw}-1)
t^{-s}
\frac{dt}{t}
+
\int_0^1
(\theta(1/t^2)^{\cw}-1)
t^{-\cw+s}
\frac{dt}{t} \\
&=
-\frac{1}{s}+\frac{1}{s-\cw} + 
\int_1^\infty
(\theta(t^2)^{\cw}-1)
t^s
\frac{dt}{t}
+
\int_0^1
(t^u\theta(t^2)^{\cw}-1)
t^{-\cw+s}
\frac{dt}{t}. 
\end{align}
Now, we observe that for $\Re(s)\ne 0$,
\begin{align}
\int_1^\infty H(-s) t^s \frac{dt}{t} &= -\frac{H(-s)}{s}\\
\int_0^1      H(s) t^s \frac{dt}{t}  &= \frac{H(s)}{s};
\end{align}
whichever integral would diverge is killed by the factor of $H(\pm s)$.
By replacing
\begin{align}
-\frac{1}{s} &= -\frac{H(s)+H(-s)}{s}\\
&=
\int_1^\infty H(-s) t^s \frac{dt}{t}
-\int_0^1      H(s) t^s \frac{dt}{t},
\end{align}
and similarly for $\frac{1}{s-\cw}$, the  formula for 
$Z_\QQ(\cw,s)$ 
simplifies to
give the desired result. ~~~$\bsq$

\noindent\paragraph{Remark.}
Since $Z_\QQ(\cw,s)$ has singularities at $s=0$ and 
$s=\Re{w}$  when $\cw =\Re (\cw)$ is real, we
cannot obtain quite as nice an expression for $\Re(w)=u$ along  vertical
lines $\Re(s) \in  \{0,\Re(\cw)\}$.
   Indeed, the Heaviside functions are precisely the contributions of
the poles as we move the integral through those points; the poles are also
reflected in the fact that for $\Re(s)\in \{0,\Re(\cw)\}$, the integral
diverges.  However, if we renormalize the integrals:
\begin{align}
\int_0^\infty f(t) \frac{dt}{t}
&\to
\lim_{B\to\infty} \frac{1}{\log B}
\int_1^B \int_{1/A}^{A} f(t) \frac{dt}{t}\, \frac{dA}{A},\\
\int_{-\infty}^\infty f(v) dv
&\to
\lim_{B\to\infty} \frac{1}{B}
\int_0^B \int_{-A}^A f(v) dv\, dA,
\end{align}
then the formula is in fact valid for all $s \neq \{ 0, u \}$; 
to prove this, use the
identity
\[
-1/s = 
\lim_{B\to\infty} \frac{1}{\log B}
\int_1^B \int_{1/A}^{A}
\frac{1}{2}\sgn(t-1) t^s \frac{dt}{t}\, \frac{dA}{A},
\]
valid for $\Re(s)=0$, and proceed as before.

%
%
%

\section{Fourier Coefficient Estimates}
\hsp
The function $\theta(t)= \sum_{n \in \ZZ} e^{-\pi n^2t}$ 
is a modular form of weight $\frac{1}{2}$
in the variable $\tau = i t $ for $\tau$ in the upper half-plane,
with a multiplier system with respect to the 
{\em theta group}\footnote{$\Gamma_\theta$ is the set of
$\left( \begin{array}{cc}
a & b \\ c & d \end{array}\right) \equiv
\left( \begin{array}{cc}
1&0 \\ 0 & 1 \end{array} \right)$
or $\left( \begin{array}{cc}
0 & 1 \\ 1 & 0 \end{array} \right)$ $(\bmod~2)$ in $PSL (2, \ZZ )$.}
$\Gamma_\theta$, 
a non-normal subgroup of index $3$
in the modular group $PSL(2, \ZZ)$ which contains $\Gamma(2)$,
the principal congruence subgroup of level $2$. 
Thus $\theta(t)^{\cw}$ is a modular form of (complex)
weight $\frac{\cw}{2}$ with (non-unitary) multiplier system
on the same group. We consider its Fourier expansion at the
cusp $i\infty$ (of width $2$) , given by
\begin{equation}\label{S301}
\theta ( it )^{\cw} = 1 + \sum_{m=1}^\infty c_m (\cw) e^{- \pi im t} ~.
\end{equation}
(The theta group has two cusps, with the
second cusp at $-1$, see Bruggeman\cite[Chap. 14]{Br94};
we do not consider the other cusp here.)
In this section our  object is to obtain estimates for the size of the
Fourier coefficients $c_m(\cw)$ as $m \to \infty$ for fixed $u$. At 
the end of the section we give explicit formulas for
a few integer values of $\cw$ where the Fourier coefficients
have arithmetic significance, namely $\cw = 0, 1, 2, 4, 6,$ and $8$.

%
%
%

\subsection{Fourier Coefficient Formulas}

We establish 
basic properties of the Fourier coefficients $c_m(\cw)$ as a
function of $\cw= u + iv$.

\begin{theorem}\label{th31}
The Fourier coefficient $c_m (\cw)$ is
a polynomial in $\QQ[\cw]$ of degree $m$.
For each $m \geq 1$, the polynomial
\begin{equation}\label{S302}
c_m^* (\cw) := (-1)^m m! ~c_m (-\cw) , \quad m \ge 1 \,,
\end{equation}
has nonnegative integer coefficients, lead term $2^m {\cw}^m$,
and vanishing constant term.
\end{theorem}

To prove this result we will need  the triple product formula
of the Jacobi theta function
$\vartheta_3 (z, q)$, see Andrews~\cite[Theorem 2.8]{A76}
or Andrews, Askey and Roy~\cite[Section 10.4]{AAR99}.

\begin{prop}~\label{pr31}
$($Jacobi Triple Product Formula$)$
The Jacobi theta function 
$$\vartheta_3 (z, q) := \sum_{n \in \ZZ} e^{2 \pi i n z} q^{n^2}$$
is given by
\begin{equation}\label{S301a}
\vartheta_3 (z, q) = \prod_{n = 1}^\infty 
(1 - q^{2n})(1 + e^{2\pi i z} q^{2n - 1})(1 + e^{-2\pi i z} q^{2n - 1}).
\end{equation}
This formula is valid for $|q| < 1$ and all $z \in \CC.$
\end{prop}

\noindent\paragraph{Proof of Theorem~\ref{th31}.}
The Fourier coefficients $c_m(\cw)$ are computable using the expansion
\begin{eqnarray}\label{S303}
\theta (t)^{\cw} & = & \left( 1 + 2 \sum_{n=1}^\infty
e^{- \pi n^2 t} \right)^{\cw} \nonumber \\
& = & 1 + \sum_{j=1}^\infty {\binom{\cw}{j}} 2^j
\left( \sum_{n=1}^\infty e^{-\pi n^2 t} \right)^j ~.
\end{eqnarray}
Terms involving $e^{- \pi mt}$ can appear only for
$1 \le j \le m$, hence we find that $c_m (\cw)$ is a polynomial 
of degree $m$ in $\cw$ with rational coefficients and
leading term $\frac{2^m}{m!} {\cw}^m$.
Clearly
\begin{equation}\label{S303a}
c_1 (\cw) = 2 \cw
\end{equation}
and
$c_2 (\cw) = 2\cw (\cw-1) ~.$
Multiplication by $m!$ clears denominators 
$m! {\binom{\cw}{j}} \in \ZZ [\cw]$ for $1 \le j \le m$ hence
$c_m^{*} (\cw)= (-1)^m m! c_m(\cw) \in \ZZ [\cw]$.

It remains to show nonnegativity.
We have
\begin{equation}\label{S304}
\vartheta_3 (0,-q)^{-\cw} = 1+ \sum_{m=1}^\infty
\frac{c_m^* (\cw)}{m!} q^m ~.
\end{equation}
The Jacobi triple product formula gives
\begin{equation}\label{S305}
\vartheta_3 (0, q) = 
\prod_{k=1}^\infty (1+q^{2k-1} )^2 (1-q^{2k} ) 
\end{equation}
hence
$$
\vartheta_4 (0, q) = \vartheta_3 (0, -q) =
(1-q) \prod_{k=1}^\infty (1-q^{2k-1} ) (1-q^{2k} ) (1-q^{2k+1} ) \,.
$$
For $\cw = u > 0$ real, we have
\begin{equation}\label{S306}
\vartheta_3 (0, -q)^{-u} = (1-q)^{-u} \prod_{k=1}^\infty
(1-q^{2k-1} )^{-u} (1-q^{2k} )^{-u} (1-q^{2k+1} )^{-u} \,.
\end{equation}
Now
\begin{eqnarray}\label{S307}
(1-z)^{-\cw} & = & 1+ \sum_{k=1}^\infty {\binom{\cw+l-1}{l}} z^l \\
\label{S308}
& = & 1+ \sum_{l=1}^\infty \sum_{m=0}^l c_{lm} {\cw}^m z^l
\end{eqnarray}
is a bivariate power series with all $c_{lm} \ge 0$.
This nonnegativity property is preserved under 
multiplication of power series, hence
$\vartheta_3 (0,-q)^{-u}$ inherits this property by \eqn{S306}.
Thus all the coefficients of $c_m^* (\cw)$ are nonnegative.
One has
\begin{equation}\label{S309}
c_m^* (\cw) = \sum_{k=0}^m c_{mk}^* {\cw}^k ~,
\end{equation}
in which $c_{m0}^* =0$ and $c_{mk}^* > 0$ for 
$1 \le k \le m$, with $c_{mm}^* = 2^m$.~~~$\bsq$

%
%
%

\subsection{Fourier coefficients for positive real $\cw$}

Let $\cw =u \in \RR_{> 0}$.
In this case $\theta (t )^u$ is a modular form of real weight $u/2$,
on the theta group $\Gamma_\theta$, with a unitary multiplier system.
Classical estimates of Petersson \cite{Pe58}  and Lehner~\cite{Le61}
for the 
Fourier coefficients of arbitrary
modular forms of positive real weight (with multiplier systems) 
show they grow polynomially in $m$, with
$$c_m (u) =O (m^{u/2 - 1}) \quad\mbox{if}\quad u>4\,,
$$
with $c_m (u) = O(m^{u/2-1} \log m )$ for $u =4$ and 
$c_m (u) = O (m^{u/4} )$ for $0 < u < 4$, where the $O$-symbol 
constants depend on $u$ in an unspecified manner.

Here we establish some weaker estimates, whose merit is
that the dependence on
$u$ is completely explicit, for use in \S6.

\begin{theorem}\label{th33}
Suppose $\cw= u \ge 0$ is real.

(i) For $m \ge 2$,
\begin{equation}\label{S313}
|c_m (u) | \le 24m^{\frac{u}{2}} ~,
\end{equation}

(ii) For $m \ge 1$,
\begin{equation}\label{S313a}
|c_m (u) | \le 6um^{\frac{u}{2}+1 } ~.
\end{equation}
\end{theorem}

\paragraph{Proof.}
(i) Write
$\vartheta_3 (q) = \vartheta_3 (0,q) = \sum_{n \in \ZZ} q^{n^2}$, so that
$\theta (t) = \vartheta_3 (e^{-\pi t} )$.
Using Cauchy's theorem we obtain the following formula:
\begin{align}
c_m(u) &=   
\frac{1}{2 \pi i} \oint \vartheta_3 (q)^u q^{-m-1} dq \nonumber \\
&=
\frac{1}{2\pi}
\int_{-\pi}^\pi \vartheta_3(R e^{i\theta})^u R^{-m} e^{-im\theta} d\theta
\nonumber \\
&
\le
R^{-m}
(\max_{-\pi\le \theta\le \pi} |\vartheta_3(R e^{i\theta})^u|),
\end{align}
for any choice $0<R<1$.  We take $R=e^{-1.01\frac{\pi}{m}}$, and thus
\[
R^{-m} = e^{1.01 \pi} < 24.
\]

For the other term, we first observe that, since $u\ge 0$,
\[
\max_{-\pi\le \theta\le \pi} |\vartheta_3(R e^{i\theta})^u|
=
(\max_{-\pi\le \theta\le \pi} |\vartheta_3(R e^{i\theta})|)^u.
\]
Since the coefficients of $\vartheta_3(q)$ are positive,
the maximum occurs for $\theta=0$; thus we must estimate $\vartheta_3(R)$.
Using the functional equation of $\vartheta_3(q)$, we find, for $m \ge 2$:
\begin{align}
\vartheta_3(R) &= 
\vartheta_3(e^{-1.01 \frac{\pi}{m}}) \nonumber \\
&= 
\sqrt{\frac{m}{1.01}}  \vartheta_3(e^{-\frac{\pi m}{1.01}}) \nonumber \\
&\le 
\sqrt{m},
\end{align}
using
\[
\vartheta_3(e^{-{\frac{2\pi}{1.01}}}) < \sqrt{1.01}.
\]
The bound \eqn{S313} follows immediately.

(ii) For $m = 1$ we have
$c_1(u)=2u \leq 6u$, so we may restrict our attention to
$m\ge 2$.
Differentiating $\vartheta_3(q)^u$ by 
$\frac{\partial}{\partial q},$ (denoted by $'$) 
and dividing by $u$, we find:
\[
q\vartheta_3'(q) \vartheta_3(q)^{u-1}
=
\sum_{m\ge 0} \frac{m c_m(u)}{u} q^m.
\]
Using Cauchy's theorem we obtain the following formula:
\begin{align}~\label{N319}
m \frac{c_m(u)}{u} &=
\frac{1}{2\pi i} \int_{|q|=R} q\vartheta_3'(q) \vartheta_3(q)^{u-1} 
\frac{dq}{q^{m+1}} \nonumber \\
& =
\frac{1}{2\pi}
\int_{-\pi}^\pi \frac{R e^{i\theta} \vartheta_3'(R e^{i\theta})}
{\vartheta_3(Re^{i\theta})} \vartheta_3(R e^{i\theta})^u R^{-m} 
e^{-im\theta} d\theta \nonumber\\
& \le
R^{-m}
(\max_{-\pi\le \theta\le \pi} |\vartheta_3(R e^{i\theta})^u|)
(\max_{-\pi\le \theta\le \pi} 
\left|\frac{R e^{i\theta} \vartheta_3'(R e^{i\theta})}
{\vartheta_3(R e^{i\theta})}\right|),
\end{align}
Taking $R=e^{-1.01\frac{\pi}{m}}$ as before, we need only estimate the third
factor.

The Jacobi triple product formula implies that the coefficients of
$q\frac{ \vartheta_3'(q)}{\vartheta_3(q)}$
are alternating in sign; in other words, using 
$\vartheta_4(q) = \vartheta_3(-q)$ the power series expansion in $q$ of
\[
-q \frac{ \vartheta_4'(q)}{\vartheta_4(q)}
\]
has positive coefficients.  In particular, its maximum is given by
\begin{align}
-R \frac{ \vartheta_4'(R)}{\vartheta_4(R)}
&= 
\frac{m^2}{(1.01)^2} e^{-\frac{\pi m}{1.01}}
\frac{\vartheta'_2(e^{-\frac{\pi m}{1.01}})}
{\vartheta_2(e^{-\frac{\pi m}{1.01}})} -
\frac{m}{2.02\pi} \nonumber \\
&\le
\frac{m^2}{(1.01)^2}e^{-\frac{\pi m}{1.01}} 
\frac{\vartheta'_2(e^{-\frac{\pi m}{1.01}})}
{\vartheta_2(e^{-\frac{\pi m}{1.01}})}.
\end{align}
Writing
\[
\frac{q \vartheta'_2(q)}{\vartheta_2(q)} = \sum_{k\ge 0} e_k q^k,
\]
we estimate
\begin{align}
e^{-\frac{\pi m}{1.01}}
\frac{\vartheta'_2(e^{-\frac{\pi m}{1.01}})}
{\vartheta_2(e^{-\frac{\pi m}{1.01}})}
&\le
\sum_{k\ge 0} |e_k| e^{-\frac{\pi k m}{1.01}} \nonumber \\
&\le
\sum_{k\ge 0} |e_k| e^{-\frac{2\pi k}{1.01}} = .25000790 \nonumber \\
\end{align}
We conclude that
\[
-R \frac{ \vartheta_4'(R)}{\vartheta_4(R)}
\le .24508175 < 1/4.
\]
Substituting  this in \eqn{N319}, we conclude
\[
m \frac{c_m(u)}{u} \le 6 m^{u/2+2}.
\]
Now \eqn{S313a} follows after multiplying  by $\frac{u}{m}$.~~~$\bsq$

Theorem~\ref{th33} implies that for real positive $\cw=u$ the
Dirichlet series
\begin{equation}~\label{S310a}
 D_u(s) = \sum_{m = 1}^\infty c_m(u) m^{-s}
\end{equation}
converges absolutely on the half-plane $\Re(s) > \frac{u}{2} + 1.$
The estimate of Petersson~\cite{Pe58}
 for the Fourier coefficients implies that  the
Dirichlet series  converges absolutely
in the half-plane $\Re(s) > \frac{u}{2}$ for $u > 4$.

It seems likely that for general $\cw \in \CC$
the Dirichlet series has no half-plane of absolute convergence,
except when $\cw=u \in \RR_{\ge 0}$, because the Fourier
coefficients grow too rapidly off the positive real axis.
We do not address this question in this paper.

%
%
%

\subsection{Maximum Size of Fourier Coefficients}

By Theorem~\ref{th31} the  maximium size of $c_m(u)$ on the
circle  $|u| = R$ occurs on the negative real axis $u = -R.$
Convergent series are known for Fourier
coefficients of modular forms of negative real integer weight,
see Petersson~\cite{Pe54} and Lehner~\cite[Theorem 1]{Le59};
the convergent series of Rademacher for the 
partition function is an example. 
The following proposition
extracts the main term in that convergent expansion,
as given in Lehner~\cite[Theorem 1]{Le59}; one 
can also prove it following the proof for the partition
function in  Apostol~\cite[Chap. 5]{Ap76}.

\begin{prop}\label{pr34}
For real $-u < 0$, there holds
\begin{equation}\label{SP320}
c_m (-u)
= (-1)^m 2 \pi u^{\frac{u}{2}}
\frac{1}{(4m)^{u+1}} I_{u+1} (\pi \sqrt{um} ) +
O (e^{(\pi + \ep ) \sqrt{max(\frac{um}{9}, (u-8)m)}} ) \,,
\end{equation}
where $I_\af (x)$ is the modified Bessel function of the
first kind, given by
\begin{equation}\label{SP321}
I_\af (x) = \left( \frac{x}{2} \right)^{\af} \sum_{k=0}^\infty
\frac{\left( \frac{x}{2} \right)^{2k}}{k!\Gamma ( \af + k+1)}.
\end{equation}
For fixed $\af > 0$,
this function satisfies 
\begin{equation}\label{SP322}
I_\af (x) = \frac{e^x}{\sqrt{2 \pi x}} \left( 1+ O_\af 
\left( \frac{1}{x} \right) \right) \quad\mbox{as}\quad
x \to \infty \,.
\end{equation}
\end{prop}

The formula  \eqn{SP320} implies a general upper bound of the form
\begin{equation}\label{SP323}
| c_m (\cw) | \le C_0 R^{\frac{R}{2}} e^{\pi \sqrt{Rm}} \,.
\end{equation}
for all complex $|\cw| = R.$

%
%
%
%

\subsection{Euler Products}

For a few  special
values of $\cw$ the associated Dirichlet series
has an Euler product, and 
the Fourier coefficients have an explicit description. 
For  $\cw \in \CC$ we define the function
$\tilde{D}_{\cw}(\frac{s}{2})$ by 
$$
\xi_{\QQ}(\cw, s) = \frac{1}{2}s(s- \cw)\pi^{-\frac{s}{2}} 
\Gamma(\frac{s}{2})\tilde{D}_{\cw}(\frac{s}{2}). 
$$
The function $\tilde{D}_{\cw}(s)$ can be assigned  a (formal)
Dirichlet series which  for $\cw  \neq 0$ is
$$
\tilde{D}_{\cw}(s) := \frac{1}{2\cw}D_{\cw}(s) \sim 
1 + \sum_{m=2}^\infty \tilde{c}_m(\cw) m^{-s}, 
$$
and for $\cw=0$ is
$$
\tilde{D}_0(s) = 1 + \sum_{m=2}^\infty \tilde{c}_m(0) m^{-s},
$$
in which 
$\tilde{c}_m(0) = \frac{c_m(\cw)}{2\cw} |_{\cw = 0}
= \frac{1}{2}\frac{d}{d\cw} c_m(\cw) |_{\cw = 0}.$
Here the use of $\sim$ indicates that the (formal) Dirichlet series expansion
based on Fourier coefficients need not have any region of
absolute convergence. However for
real nonnegative $\cw$  it does converge on a half-plane,
as follows from Theorem~\ref{th33}. 

It is easy to determine which values $\cw$ give (formal) Dirichlet series
that have an Euler product.

\begin{lemma}~\label{lem35}
For $\cw \in \CC$ the  formal Dirichlet series 
assigned to $\tilde{D}_{\cw}(s)$ has an Euler
product expansion if and only if $\cw = 0, 1, 2, 4$ and $8$.
\end{lemma}

\noindent\paragraph{Proof.}
To see that $\cw = 0, 1, 2, 4,$ and $8$ are the only
complex values for which $\tilde{D}_{\cw}(s)$ can have
an Euler product,
we consider the necessary conditon
$$
\tc_2(\cw)\tc_3(\cw) = \tc_6(\cw).
$$
This gives a polynomial of degree five in $\cw$ whose roots are
$\cw = 0, 1, 2, 4,$ and $8$. 

For $\cw = 1, 2, 4$ and $8$ the Dirichlet series $\tilde{D}$
has an Euler product. For $\cw = 1$ we have already seen that
$$ \tilde{D}_1(s)  = \zeta(2s) =  \prod_{p}(1 - p^{-2s})^{-1}. $$ 
For $\cw = 2$, 
$$ 
\tilde{D}_2(s) =  \sum_{m= 1}^\infty  (-4 / m) m^{-s}    
         = (1 - 2^{-s})^{-1}\prod_{p=1(4)} (1-p^{-s})^{-1}
           \prod_{p=3(4)} (1+p^{-s})^{-1}.
$$
For $\cw = 4,$
\begin{eqnarray}
\tilde{D}_4(s) & = &
(1-4 \cdot 2^{-2s}) \zeta(s) \zeta(s-1) \nonumber\\
& = &
 \frac{1+2\cdot  2^{-s}}{1-2^{-s}}
      \prod_{p ~odd} (1-p^{-s})^{-1}(1-p \cdot p^{-s})^{-1}.\nonumber
\end{eqnarray}
For $\cw = 8$,
\begin{eqnarray}
\tilde{D}_8(s) & = &  (1-2\cdot 2^{-s}+16\cdot 2^{-2s}) \zeta(s)\zeta(s-3) 
\nonumber\\
 & = &  \frac{1-2 \cdot  2^{-s}+16\cdot  2^{-2s}} 
{(1-2^{-s}) (1-8 \cdot 2^{-s})} 
\prod_{p ~odd} (1-p^{-s})^{-1}(1-p^3\cdot p^{-s})^{-1}. \nonumber
\end{eqnarray}
These Dirichlet series  are scaled multiples of the zeta functions for
the rational field, the Gaussian integers, the integral quaternions
and the integral octaves, respectively. 

For
 $\cw = 0$ the Dirichlet series also has an Euler product, which is
\begin{equation}~\label{eq324}
 \tilde{D}_0(s)=
\frac{1 -  2^{1-s}}{1 - 2^{- s}}
\prod_{p~odd}(1 - p^{-s})^{-1}(1 - p^{-1 - s})^{-1}.
\end{equation}
This follows from Theorem~\ref{nth51} below; note that it converges
absolutely for $\Re(s) > 1.$
~~~$\bsq$

One immediately sees from the expression for  $\xi_{\QQ}(4, s)$ 
as a product of shifted Riemann zeta functions that most of
its zeros cannot be on  ``critical line'' $\Re(s) = 2.$ The only zeros
on the line come from the Euler factor at the prime $2$ and have
are of number  $O(T)$ up to height $T$, while there are $\Omega(T\log T)$
zeros off the line coming from the Riemann zeta zeros. Exactly the
same thing happens for $\xi_{\QQ}(8, s).$ 

There are other positive integer values of $\cw$ where the Dirichlet series
can be determined explicitly, without an Euler product.
For $\cw = 6$ the Dirichlet series is a linear combination of two
Dirichlet series with Euler products, namely 
$$
\tilde{D}_6(s)= \frac{4}{3}\zeta(s-2)L(s, \chi_{-4}) - 
\frac{1}{3} \zeta(s)L(s-2, \chi_{-4}),
$$
in which 
$L(s, \chi_{-4}) = \sum_{m= 1}^\infty  (-4 / m) m^{-s},$
with $(-4/ m)$ being the Jacobi symbol.

%
%
%

\section{Growth Bounds}

%
%
%

It is well known that the Riemann $\xi$-function
$\xi(s) = \frac{1}{2} s (s - 1) \pi^{\frac{s}{2}} \Gamma(\frac{s}{2})\zeta(s)$
is an entire function of order one and infinite type. It
is bounded in vertical strips. In this section we show that both
these properties generalize to the function $\xi_{\QQ}(\cw,s)$.

\subsection{Growth of Maximum Modulus}

We prove the following bound for the  the two-variable
zeta function.

\begin{theorem}~\label{Nth41}
There is a constant $C>0$ such that the 
entire function $\xi_\QQ(\cw,s)$ 
satisfies the
growth bound: For all $(\cw, s) \in \CC^2$, if  $R = |\cw|+|s|+1$ then
\begin{equation}\label{S212}
|\xi_\QQ(\cw,s)| \le e^{C R \log R}.
\end{equation}
\end{theorem}

\noindent\paragraph{Proof.}
Set
\[
I(\cw, s) = \int_0^1 \frac{\theta(\frac{1}{t^2})^{\cw}-1}{\cw}t^{-s} 
\frac{dt}{t}
= \frac{1}{2} \int_1^\infty \frac{\theta(t)^{\cw}-1}{\cw} t^{s/2} \frac{dt}{t}
\]
Then (2.9) gives
\[
2\xi_\QQ(\cw,s) = 1 + s(s-\cw)(I(\cw, s)+I(\cw, \cw-s)).
\]

For $t > 0,$
\[
\theta(t)
=
1 + 2\sum_{m\ge 1} e^{-\pi m^2 t},
\]
and $\theta(t)>1$ for all $t>0$
since all of the terms are positive.
It follows that for all positive $t$ and all $ \cw \in \CC$,
\[
\left|\frac{\theta(t)^{\cw}-1}{\cw}\right|
\le
\frac{\theta(t)^{|\cw|}-1}{|\cw|}.
\]
Moreover, the right-hand-side is increasing in $|\cw|$, so we find
that for all $t > 0,$
\[
\left|\frac{\theta(t)^{\cw}-1}{\cw}\right|
\le
\frac{\theta(t)^R - 1}{R}
\le \theta(t)^R   \log\theta(t),
\]
using 
$\frac{e^{Rx}-1}{R} \le x e^{Rx}$
for $x\ge 0$. 

Similarly, the function $e^{\pi t}(\theta(t)-1)$ is
decreasing in $t$, so we find that for $1\le t<\infty$, we have:
\[
\theta(t) \le 1 + C_1 e^{-\pi t} \le e^{C_1 e^{-\pi t}},
\]
where $ C_1 = e^\pi (\theta(1)-1)$.
In particular, for $t \geq 1$,
\[
 \theta(t)^R \log\theta(t)
\le
e^{C_1 e^{-\pi t} R} C_1 e^{-\pi t}
=
C_1 e^{C_2 R} e^{-\pi t}.
\]

We thus have
\begin{align}
|I(\cw, s)|
&\le
\frac{1}{2} \int_1^\infty
\left|\frac{\theta(t)^{\cw}-1}{\cw}\right| |t^s| \frac{dt}{t} \nonumber\\
&\le
\frac{1}{2} 
C_1 e^{C_2 R}
\int_1^\infty
e^{-\pi t} t^{\frac{|s|}{2}} \frac{dt}{t} \nonumber \\
&\le
\frac{1}{2} 
C_1 e^{C_2 R}
\int_0^\infty
e^{-\pi t} t^{\frac{R}{2}} \frac{dt}{t} \nonumber \\
&=
\frac{1}{2} C_1 e^{C_2 R} \pi^{-\frac{R}{2}} \Gamma(\frac{R}{2}).
\end{align}
Since $R\ge 1$, and $\Gamma(R) \le R^R$ on this
range, it follows that there exists
a constant $C_3>0$ such that
\[
|I(\cw, s)|\le e^{C_3 R\log R}.
\]
Since both $1$ and $|s(s - \cw)|$ are $o(R^R)$, the theorem follows.~~~$\bsq$

A notion of entire function of finite order for
functions of several complex variables is described in
Ronkin \cite{Ro74}, \cite{Ro89} and Stoll~\cite{St74}.
In particular, there is a Weierstrass
factorization theorem for such functions in terms of their
zero locus.

The zero locus $\sZ_\QQ$ of $\xi_\QQ (\cw,s)$ is a 
one-dimensional complex analytic
manifold, possibly with singular points, which may have many connected
components. One way to study it is to take ``linear slices'' 
to obtain functions
$h(s)$ of one complex variable,
whose zero sets consist of  isolated points which can be counted.
If $\alpha := (\af_1, \af_2 )$ and $\beta := (\beta_1, \beta_2 )$
with $\alpha_j, \beta_j \in \CC$, the
{\em linear slice function}
$h_{\alpha , \beta} (s)$ of $\xi_\QQ (\cw,s)$, is
\begin{equation}\label{S401}
h_{\af, \beta} (s) := \xi_\QQ (\af_1 s + \af_2, \beta_1 s + \beta_2 ) ~,
\end{equation}
where we assume $| \af_1 | + | \beta_1 | \neq 0$ to avoid constant
functions.

\begin{lemma}\label{le41}
Any linear slice function $h_{\alpha, \beta} (s)$ is an 
entire function of order at most $1$.
If $N_{\alpha, \beta}^\ast (T) = \# \{\mbox{zeros of $h_{\alpha, \beta} (s)$ 
with $|s| \le T$} \}$, 
then
\begin{equation}\label{S402}
N_{\af, \beta}^\ast (T) = O(T \log T ) \quad\mbox{for}\quad T \ge 2 ~,
\end{equation}
where the implied constant in the $O$-symbol depends on $\alpha, \beta$.
\end{lemma}

\paragraph{Proof.}
This follows from the growth estimate in 
Theorem \ref{Nth41}, using
Jensen's formula.~~~$\bsq$

This result applies in particular to  linear slices where
$\cw \in \CC$ is held fixed, i.e. 
$\af_1 =0$, $\af_2 = \cw$, 
$\beta_1 =1$ and $\beta_2 =0$.
 which gives the function $ \xi_{\QQ} (\cw,s),$
with $\cw$ regarded as constant, e. g. in \S7.

%
%
%

\subsection{Growth Bounds on Vertical Lines}

We next consider growth bounds for $\xiq(\cw, s)$ on vertical
lines $s = \sigma +it$ with $\sigma$ held fixed.
Recall that a function $f(x)$ is in the {\em Schwartz space} 
$\sS(\RR)$ if and only if for each
$m, n  \geq 0$ there is a finite constant $C_{m, n}$ such that
$$ \sup_{x \in \RR}   |x^n \frac{d^m}{dx^m} f(x)| \leq C_{m, n},$$
with a similar definition for functions defined on a closed half-line
$\RR_{\geq 0}$ or $\RR_{\leq 0}$.

\begin{theorem}~\label{th43}
For each $\cw \in \CC$, $\sigma \in \RR$, and real
$-\frac{\pi}{4}< y < \frac{\pi}{4}$, the function
\[
e^{yt}\xi_{\QQ}(\cw, \sigma + it) ~~~~~~  -\infty < t < \infty,
\]
belongs to the Schwartz space $\sS(\RR).$ Furthermore,
the implied Schwartz constants can be chosen
uniformly on any compact subset  
$\Omega $ of 
$ \{( \cw, \sigma, y): \cw \in \CC, \sigma \in \RR,  
-\frac{\pi}{4}< y < \frac{\pi}{4} \}.$
In particular, these functions are bounded in vertical strips, i.e.
there is a finite constant $A_1 = A_1(\Omega)$ such that
$$
\sup_{(\cw, \sigma, y) \in \Omega,~ t \in \RR} 
|e^{yt}\xi_{\QQ}(\cw, \sigma + it)| \leq A_1.
$$
\end{theorem}

We will deduce this result from an integral representation 
of  $\xi_{\QQ}(\cw, \sigma + it),$ which we prove first. 
We define the function $h(\cw, z)$  
for
$\cw \in \CC$ and $z = x + iy \in \CC$ with 
$-\frac{\pi}{4} < y < \frac{\pi}{4},$ by
\[
h(\cw, z) := 
\begin{cases}
\frac{\theta(e^{2z})^{\cw} -1}{\cw} & \cw \ne 0\\
\log\theta(e^{2z}) & \cw =0.
\end{cases}
\]
Next we define
\begin{equation}~\label{eq44a}
\gaa(\cw, z)= ( \frac{d^2}{dz^2} + w \frac{d}{dz}) h(w, z).
\end{equation}
We have the following Fourier-Laplace transform formula for
$\xi_{\QQ}(\cw, s),$ valid on $\CC^2.$

\begin{theorem}~\label{th44}
Let $\cw \in \CC$, and real $y$ with
$-\frac{\pi}{4} < y < \frac{\pi}{4}.$ Then for all
$s \in \CC$,
\[
 \xi_{\QQ}(\cw, s) =  \frac{1}{2}\int_{-\infty}^\infty 
 \gaa(\cw, x + iy) e^{xs} dx.
\]
Regarded as a Fourier transform, with $s = \sigma + it$
with fixed $\sigma \in \RR$, the  integrand 
\[
f_{w, \sigma, y}(x) := 
e^{\sigma x} \gaa(\cw, x + iy), ~~~~~
 -\infty < x < \infty,
\]
belongs to the Schwartz space $\sS(\RR)$ and the implied Schwartz
constants are uniform on compact subsets of
$\{ (w, \sigma, y): \cw \in \CC, \sigma \in \RR, 
-\frac{\pi}{4} < y < \frac{\pi}{4} \}.$

\end{theorem}

To prove Theorem~\ref{th44} we formulate several estimates as 
preliminary lemmas.

\begin{lemma}\label{lem45}
Define an analytic function $f(\cw, \tau)$ for $\cw, \tau  \in \CC$, 
with $\Re(\tau)>0$, by
\[
f(\cw, \tau) = 
\begin{cases}
\frac{\theta(\tau)^{\cw} -1}{\cw} & \cw \ne 0\\
\log\theta(\tau) & \cw =0.
\end{cases}
\]
For any integer $k\ge 0$, 
and for any region of the form 
$$\Omega_{R, \epsilon} = \{ (\cw, \tau):~ 
|\cw|\le R ~\mbox{and}~\Re(\tau)\ge \epsilon\},$$
there is a finite constant $C_k = C_k(R, \epsilon)$ such that
\[
|e^{\pi \tau} \frac{d^k}{d{\tau}^k} f(\cw,\tau)| \le C_k.
\]
\end{lemma}

\noindent\paragraph{Proof.}
For $\cw \neq 0$ the results of  \S3 show
 that $f(\cw, \tau)$ has a Fourier expansion:
\[
f(\cw, \tau) = \sum_{m= 1}^\infty \frac{c_m(\cw)}{\cw} e^{-m\pi \tau}.
\]
Note that $\frac{c_m(\cw)}{\cw}$ is well-defined as a polynomial in $w$,
hence this expansion makes sense for $\cw = 0$ as well,
with Fourier coefficient  $\frac{d}{d\cw}c_m(\cw)|_{\cw = 0}.$
Now, for $m\ge 1$, we have the inequality
\[
\left|\frac{c_m(\cw)}{\cw}\right| \le \frac{(-1)^m c_m(-|\cw|)}{|\cw|};
\]
this follows from the fact that $(-1)^m c_m(-\cw)$ has positive coefficients,
and is 0 at $\cw=0$.  Moreover, this upper bound is monotonically increasing
in $|\cw|$.  
In particular, since $\vartheta_3(q)$ has radius of convergence 1,
and no zeros in the unit disk, we find that for fixed  $|\cw|>0$,
\[
\frac{\vartheta_3(q)^{-|\cw|}-1}{|\cw|}
\]
has radius of convergence 1, and thus
\[
\lim_{m\to\infty} \left|\frac{c_m(\cw)}{\cw}\right|^{1/m} = 1.
\]
In particular, since $q = e^{- \pi \tau}$, the condition
$\Re(\tau) > 0$ makes $f(\cw, \tau)$ analytic in the region 
$\Omega_{R, \epsilon}.$

Similarly, for the $k$th derivative, we have:
\[
\frac{d^k}{dt^k} f(\cw, \tau)
=
\sum_{m=1}^\infty (-\pi m)^k \frac{c_m(\cw)}{\cw} e^{-m\pi \tau},
\]
and, using term-by-term absolute value estimates,
\[
|e^{\pi \tau} \frac{d^k}{d{\tau}^k}f(\cw, \tau)| \leq
\sum_{m=1}^\infty (\pi m)^k \frac{(-1)^m c_m(-|\cw|)}{|\cw|} 
e^{-(m-1)\pi \Re(\tau)}.
\]
Over the region $|\cw|\le R$, $\Re(\tau) \ge \epsilon$, this is
bounded by its value for $\cw=R$, $\tau = \epsilon$; as the sum converges
for those values, we obtain the required uniform bound. 
~~~$\bsq$

\begin{lemma}\label{lem46}
Let $z= x + iy$, and for
fixed  $-\frac{\pi}{4} < y < \frac{\pi}{4}$ 
and fixed $\sigma \in \RR$, the function 
\[
f_{\cw,\sigma,y}(x) = e^{\sigma (x + iy)} \gaa(\cw, x + iy)~~~~~ 
 -\infty < x < \infty,
\]
where $\gaa(w, z)$ is defined as
in Theorem~\ref{th44}, 
lies in the Schwartz space $\sS(\RR)$. Moreover, the
implied Schwartz bounds are uniform over compact regions in 
$\{ (\cw, \sigma, y) : \cw \in \CC, \sigma \in \RR, 
 -\frac{\pi}{4} < y < \frac{\pi}{4}\}$.
\end{lemma}

\noindent\paragraph{Proof.}
By definition 
$h(\cw, z) = f(\cw, e^{2z})$,
in which  $-\frac{\pi}{4} < y= \Im(z)  < \frac{\pi}{4}$ so that
$\Re(\tau) = \Re{e^{2z}} > 0.$ We claim that for any $\sigma \in \RR$,
the function
$\{ e^{\sigma (x + iy)}h(\cw, x+iy)~:~ 0 \leq x < \infty\}$ 
is a Schwartz function on the half-line $\sS(\RR_{\geq 0})$, and that the
implied Schwartz bounds are uniform over compact regions in
$\cw \in \CC$, $-\frac{\pi}{4} < y  < \frac{\pi}{4}$, $\sigma \in \RR$.
To see this, note that for 
$\Re(z)\ge 0$, we find $\Re(e^{2z})\ge \cos(2\Im(z)),$
and since $\Im(z)$ is bounded away from $\pm \pi/4$,
$\Re(e^{2z})$ is bounded away from 0. 
Now using  Lemma~\ref{lem45}, we have
\[
|\Re(z)^n \frac{d^m}{dz^m} e^{\sigma z} f(\cw, e^{2z})|
=
\left||\log(\tau)|^n (\tau \frac{d}{d\tau})^m \tau^{\sigma} f(\cw, \tau^2)
\right|_{\tau=e^z}
=
O(e^{(\sigma+m)z} z^n e^{-\pi e^z})
=
o(1),
\]
as $x = \Re(z) \to \infty$, 
uniformly in a compact region 
in $\cw$, $-\frac{\pi}{4} < \Im(z)  < \frac{\pi}{4}$, $\sigma$, 
proving the claim. (The function $|f(w, e^{2z})|$ and its derivatives
go to zero at a super-exponential rate as $x \to \infty$, all other
variables fixed.)

Now applying the operator $\frac{d}{dz}(\frac{d}{dz}+w)$
to $h(w, z)$, 
the claim implies that $\gaa(\cw, x+iy)$
is a Schwartz function in $x$ on $\RR_{\geq 0}.$
The functional equation
$\theta(e^{2z}) = e^{-z} \theta(e^{-2z})$ then yields
\[
\gaa(\cw, -z)
=
e^{\cw z} \gaa(\cw,z),
\]
so that  $\gaa(\cw, x + iy)$ is a Schwartz function in $x$ on  $\RR_{\leq 0}$
as well. Thus it is in $\sS(\RR)$, and the uniformity of the estimates in
compact regions in $w$, 
$ -\frac{\pi}{4} < \Im(z) < \frac{\pi}{4}$
is inherited.
~~~$\bsq$

To prove Theorem~\ref{th44}, we will use repeated integration by
parts starting from the integral representation for 
$Z_{\QQ}(\cw, s)$ in Theorem~\ref{th22a}. To justify this step
we show that the integrand of that representation is a Schwarz
function for $\Re(s)  \notin \{ 0, \Re(\cw) \}.$

\begin{lemma}\label{lem47}
Let $\cw, z\in \CC$, and $\sigma \in \RR$ be fixed, with
$\sigma \notin \{ 0, \Re(\cw) \}$. Then for
$z = x + iy,$ and $-\frac{\pi}{4}< y < \frac{\pi}{4}$, the functions
of $x \in \RR$, 
\[
g_{\sigma}(\cw, x+ iy) =
\begin{cases}
\frac{e^{\sigma(x + iy)}}{\cw} ( \theta(e^{2(x + iy)})^{\cw} - H(\sigma)
-H(\cw - \sigma)e^{-w(x + iy)}) & \cw \ne 0  \\
e^{\sigma(x + iy)}(\log \theta(e^{2(x + iy)}) - H(\sigma) - H(-\sigma)(x + iy))
& \cw = 0
\end{cases}
\]
all belong to  the Schwartz space $\sS(\RR).$  
\end{lemma}

\noindent\paragraph{Proof.}
It suffices to show that $g_{\sigma}(\cw,x + iy)$ and
$g_{\sigma}(\cw, -x - iy)$ are Schwartz functions on the half-line
$x \ge 0$. 
We have
\[
g_{\sigma}(\cw,x + iy)
=
e^{\sigma (x+iy)} f(\cw,e^{2(x+iy)}) + \frac{1}{\cw}(H(-\sigma) e^{\sigma (x + iy)}-
H(\cw-\sigma)e^{(\sigma-\cw)(x + iy})
\]
using the relation $1 - H(\sigma)= H(-\sigma).$ We also have
\[
g_{\sigma}(\cw, -x - iy)
=
e^{(\cw-\sigma)(x + iy)} f(\cw, e^{2(x + iy)}) +
\frac{1}{\cw}(-H(\sigma) e^{-\sigma (x + iy)}+
H(\sigma-\cw)e^{(\cw-\sigma)(x + iy)}),
\]
using the functional equation for $\theta(e^{2z}).$
In each case, the first term has been shown to be Schwartz, 
and to be uniform in the parameters, in the proof of Lemma~\ref{lem46}.
It remains only to consider the ``correction'' terms, on
the half-line $x \ge 0$.

Suppose $\cw\ne 0$, it suffices to observe that 
for all $s \in \CC$ with $\Re(s) \neq 0$ and fixed $y$ the function 
$H(-s) e^{s (x + iy)}$
is Schwartz on the half-line  $x = \ge 0$;
indeed,
\[
x^n \frac{d^m}{dx^m} H(-s) e^{s(x+iy)}
=
H(-s) s^m x^n e^{s(x+iy)}.
\]
When $\Re(s)<0$, the exponential dominates, and thus the function is
bounded; when $\Re(s)>0$, $H(-s)=0$, so the function is 0.  Since $\cw\ne 0$,
we can safely divide by $\cw$ without affecting boundedness.

Suppose $\cw=0$. Then  the only new  ``correction'' term is
 the function $H(-s) (x + iy) e^{s (x + iy)}$, 
which again on a half-line $x \ge 0$ is in the  Schwartz space.
~~~$\bsq$

\noindent\paragraph{Remark.} With some further work
it can be shown that in Lemma~\ref{lem47}
the implied constants for the Schwartz functions
are uniform over compact regions in $(\cw, \sigma, y)$-space that avoid
the two lines $\sigma= 0$ and $\sigma = \Re(\cw).$ However
we do not need this result.

\noindent\paragraph{Proof of Theorem~\ref{th44}.}
We take $s=\sigma+it$.  When $\sigma\notin\{0,\Re(w)\}$, $w\ne 0$,
Theorem~\ref{th22a} gives 
\[
\xi_\QQ(w,\sigma+it)
=
\frac{1}{2}(\sigma+it)(\sigma+it - \cw)
\int_{-\infty}^\infty e^{it (x + iy)} g_{\sigma}(w,x+iy) dx,
\]
where $g_{\sigma}(w,x+iy)$ is given in Lemma~~\ref{lem47}.
Now $ g_{\sigma}(\cw,x+iy)$ is a Schwartz
function in $x$ by Lemma~\ref{lem47}, so we may integrate by parts twice,
using $z = x + iy$, obtaining
\begin{align}
\xi_\QQ(\cw,\sigma+it)
&=
\frac{1}{2}\int_{-\infty}^\infty e^{it z}
\left(\sigma-\frac{d}{dz}\right)\left(\sigma- \cw- \frac{d}{dz}\right)
g_\sigma( \cw, x + iy)dx\\
&=
\frac{1}{2}\int_{-\infty}^\infty
e^{(\sigma+it)z}
\left(-\frac{d}{dz}\right)\left(-\cw-\frac{d}{dz}\right)
f(\cw, e^{2z}) dx   \\
&=
\frac{1}{2}\int_{-\infty}^\infty e^{(\sigma+it)(x + iy)} \gaa(\cw, x+ iy) dx \\
&=
\frac{1}{2} e^{-ty + i \sigma y} \int_{-\infty}^\infty
(e^{\sigma x} \gaa(\cw, x+ iy))e^{i t x} dx.
\end{align}
This integral agrees with $\xi_\QQ(w,s)$ off the lines $\sigma = 0, \Re(\cw).$
Since the integrand is uniformly Schwartz
by Lemma~\ref{lem46}, this integral gives an analytic
function of $w$ and $s$, and must therefore agree with $\xi_\QQ(w,s)$
everywhere.
~~~$\bsq$

\noindent\paragraph{Proof of Theorem~\ref{th43}.}
View the  integral representation of $\xi_{\QQ}(w, \sigma - it)$
in
Theorem~\ref{th44} as a  
Fourier transform, with $s = \sigma + it$ with fixed $\sigma \in \RR.$
Since the Fourier transform maps Schwartz 
space $\sS(\RR)$ to itself, it follows that for
fixed $-\frac{\pi}{4} < y < \frac{\pi}{4}$, 
\[
 e^{-yt}\xi_{\QQ}(\cw, \sigma - it) ~~~~~~  -\infty < t < \infty,
\]
belongs to $\sS(\RR)$. The uniformity of the Schwartz constants
on compact subsets $\Omega$ of  $ (\cw, \sigma, y)$-space  is
inherited from the corresponding uniformity property in Theorem~\ref{th44}.
~~~$\bsq$

We conclude this section with another consequence of the Fourier-Laplace
integral representation of $Z_{\QQ}(\cw, s)$
 by (uniform) Schwartz functions.

\begin{lemma}~\label{lem48}
Let $Q(T) \in \CC[T]$ be any polynomial. Then for any $s = \sigma + it$ with 
$\sigma \notin \{0, \Re(w)\}$ and $z = x + iy$ with
$-\frac{\pi}{4}< y < \frac{\pi}{4}$,
\[
\frac{1}{2\pi}\int_{-\infty}^{\infty}
\left(Q(\sigma + it) Z_{\QQ}(\cw, \sigma + it)\right)e^{-(\sigma+ it)z}dt =
Q(-\frac{d}{dz})(\theta(e^{2z})^\cw - H(\sigma)- 
H(\cw - \sigma)e^{-\cw z})|_{z= x+iy}.
\]
Here the integrand is a Schwartz function of $t$.
\end{lemma}

\noindent\paragraph{Proof.}
Using  Lemma~\ref{lem47}, the integral is 
a Fourier transform with integrand in $\sS(\RR)$, viewing $\sigma$
as fixed. The case $ Q(z) \equiv 1$ follows from 
Theorem \ref{th22a} by taking the Fourier transform, since the
right side of that formula can be viewed as an inverse Fourier transform.
Now use the fact that the Fourier transform leaves $\sS(\RR)$
invariant and transforms multiplication by $s$ to differentiation, and
apply $Q(-\frac{d}{dz})$ to both sides of the identity
with polynomial $Q(z) \equiv 1.$   
~~~$\bsq$

%
%
%

\section{Case $\cw = 0$}

We evaluate
the entire function $\txi_{\QQ}(\cw, s)$ in
the plane $\cw = 0$.

\begin{theorem}\label{nth51}
The entire function $\txi_{\QQ}(\cw, s)$ 
of two complex variables has
\begin{equation}~\label{neq51}
\txi_{\QQ}(0, s)= -\frac{s^2}{8}(1 - 2^{1 + \frac{s}{2}})
(1 - 2^{1 - \frac{s}{2}}) \hat{\zeta}(\frac{s}{2}) 
\hat{\zeta}(\frac{-s}{2}),
\end{equation}
in which $\hat{\zeta}(s)=  \pi^{- s/2} \Gamma (s/2) \zeta (s).$
\end{theorem}

\noindent\paragraph{Remarks.}
(1) It is evident from \eqn{neq51}
that $\txi_{\QQ}(0, s)$ satisfies the functional equation
$$
\txi_{\QQ}(0, s)= \txi_{\QQ}(0, -s).
$$

(2) Comparing  the formula
$\txi_{\QQ}(0, s)= \frac{s^2}{2}\pi^{-\frac{s}{2}} \Gamma(\frac{s}{2}) 
\tilde{D}_0(\frac{s}{2})$
with Theorem~\ref{nth51},
and  using the functional equation for $\zeta(s)$ leads to 
$$
\tilde{D}_0(\frac{s}{2})= \frac{1}{4}(1 - 2^{1 + \frac{s}{2}})(1 - 2^{1 - \frac{s}{2}})
\zeta(\frac{s}{2})\zeta(1 + \frac{s}{2}).
$$
It is evident that this
 Dirichlet series has an Euler product, already stated in \S3.4. 

The proof of Theorem~\ref{nth51}
depends on the Jacobi triple product formula, which is
used in evaluating the Fourier coefficients of $\log \theta(t)$.
We state this as a preliminary lemma.

\begin{lemma}\label{nle51}
The coefficients $c_m^{'}$ of
$$ \log \theta(t) = \sum_{m=1}^\infty c_m^{'} e^{- \pi m t}$$
are given by
\begin{equation}~\label{neq52}
 c_m^{'} = 2 \sigma_{-1}(m) - 5\sigma_{-1}(\frac{m}{2}) + 
2\sigma_{-1}(\frac{m}{4}),
\end{equation}
where
\[
\sigma_{-1}(m) =
 \begin{cases}
\sum_{d | m} \frac{1}{d} & m \in \ZZ_{>0} \\
0  & \mbox{otherwise}.
\end{cases} 
\]

\end{lemma}

\noindent\paragraph{Proof.}
We use the eta product
\[
\tilde\eta(t) := \prod_{k = 1}^\infty (1-e^{-k\pi t}).
\]
By the Jacobi triple product formula, we have
\begin{align}
\theta(t)
&=
\prod_{k=1}^\infty (1+e^{-(2k-1)\pi t})^2(1-q^{-2k\pi t}) \nonumber\\
&=
\prod_{k=1}^\infty (1-e^{-k\pi t})^{-2}
\prod_{k=1}^\infty (1-e^{-2k\pi t})^5
\prod_{k=1}^\infty (1-e^{-4k\pi t})^{-2} \nonumber \\
&=
\tilde\eta(t)^{-2}\tilde\eta(2t)^5\tilde\eta(4t)^{-2}.
\end{align}
If we define  $\kappa_m$ by
$$
\log\tilde\eta(t) := \sum_{m=1}^\infty \kappa_m e^{-m\pi t},
$$
then we have
\begin{equation}~\label{neq503}
c_m^{'} = -2\kappa_m+5\kappa_{m/2}-2\kappa_{m/4},
\end{equation}
with the convention that $\kappa_m=0$ if $m\notin\ZZ$.  Now,
\begin{align}
\log\tilde\eta(t)
&=
\sum_{k=1}^\infty \log(1-e^{-k\pi t}) \nonumber \\
&=
-\sum_{k=1}^\infty \sum_{m=1}^\infty  \frac{1}{m} e^{-mk\pi t} \nonumber\\
&= 
-\sum_{m=1}^\infty (\sum_{d|n} \frac{1}{d}) e^{-n\pi t},
\end{align}
so that
$\kappa_m = - \sigma_{-1}(m)$
as required.
~~~$\bsq$

\noindent\paragraph{Proof of Theorem~\ref{nth51}.}
Since $\txi_{\QQ}(\cw, s)$ is an entire function of two variables
we have, for positive real $\cw= u$,
\[
\txi_\QQ(u,s) = \lim_{u \to 0^+} \frac{s(s-u)}{2u} Z_{\QQ}(u, s).
\]
Fix $\Re(s)>0$; for $\Re(s)> u$, we have
\[
\txi_\QQ(u,s)
=
\frac{1}{2}s(s-u) \int_0^\infty \frac{\theta(t^2)^u-1}{u} t^s \frac{dt}{t}
= 
\frac{1}{2}s(s-u) \int_0^\infty \frac{e^{u\log \theta(t^2)}-1}{u} t^s \frac{dt}{t}
\]
Suppose $\Re(s) > 0.$
Letting $u \to 0$, and using  
$\lim_{u \to 0^+} \frac{e^{\alpha u} -1}{u} = \alpha,$
we eventually have $\Re(s) > u$ and so  
we legitimately obtain
\[
\txi_\QQ(0,s)
=
\frac{s^2}{2}
\int_0^\infty
\log(\theta(t^2)) t^s \frac{dt}{t}.
\]
Now expand $\log \theta(t^2)$ in Fourier series and integrate
term-by-term to obtain
\[
\txi_\QQ(0,s)
=
\frac{s^2}{4} \pi^{-\frac{s}{2}}\Gamma(\frac{s}{2}) 
( \sum_{m=1}^\infty \ttc_m m^{-\frac{s}{2}}) 
\]
and 
this converges for $\Re(s) > 2$
since $\ttc_m = O (\log m)$
by Lemma~\ref{nle51}.
Setting 
\[
K(s) = \sum_{m=1}^\infty \kappa_m m^{-s},
\]
and then using \eqn{neq503} gives
\[
\txi_\QQ(0,s)
=
-\frac{s^2}{4} \pi^{-\frac{s}{2}}\Gamma(\frac{s}{2})
(2 - 5 \cdot 2^{-\frac{s}{2}} + 2\cdot 4^{-\frac{s}{2}}) K(\frac{s}{2}).
\]
However we have 
\begin{align}
K(s)
&=
-\sum_{m=1}^\infty (\sum_{d|m} \frac{1}{d}) m^{-s} \nonumber \\
&=
-\prod_p
\sum_{k= 0}^\infty \frac{1-p^{-k-1}}{1-p^{-1}} p^{-sk} \nonumber\\
&=
-\prod_p (1-p^{-s})^{-1}(1-p^{-s-1})^{-1}  \nonumber\\
&=
-\zeta(s)\zeta(1+s),
\end{align}
valid whenever $\Re(s)>1$ (since $\kappa_m=O(\log m)$).
We thus conclude that for $\Re(s) > 2$
\[
\txi_\QQ(0,s)
=
\frac{s^2}{4} (2 - 5\cdot 2^{-s/2} + 2 \cdot 4^{-s/2})
\pi^{-\frac{s}{2}} \Gamma(\frac{s}{2})
\zeta(\frac{s}{2}) \zeta(1+\frac{s}{2})
\]
Using the duplication formula 
$\Gamma(\frac{s}{2})= \frac{1}{\sqrt{2\pi}}2^{\frac{s}{2} - \frac{1}{2}}
\Gamma(\frac{s}{4})\Gamma(\frac{s}{4} + \frac{1}{2})$
and the functional equation for the zeta function, we obtain
for $\Re(s) >2$ that
\[
\txi_\QQ(0,s)
=
-\frac{s^2}{8}
(1-2^{1+s/2})(1-2^{1-s/2})
\hat\zeta(s/2) \hat\zeta(-s/2).
\]
Since both sides are analytic in $s$, the formula is
valid for all $s\in \CC$.
~~~$\bsq$

%
%
%

\section{Location of Zeros: $\cw$ Positive Real}
\hsp
In this section we suppose $\cw= u >0 $ is fixed, and study
the zeros of $\xi_{\QQ}(u, s)$. 
We will first show that the zeros of $f_u(s)$ are localized
in a vertical strip centered on the
``critical line'' ${\rm Re}{s} = \frac {u}{2},$ whose width depends on $u$,
and then we shall derive an estimate  for the number of zeros with
imaginary part of height at most $T$. 

%
%
%

\subsection{Vertical Strip Bound}

We show for real $\cw = u > 0$ that the zeros of $\xi_{\QQ}(u, s)$
are confined to a vertical strip of width $u + 16$. 
Theorem~\ref{nth51} implies that this bound is valid for $u = 0$ as well.

\begin{lemma}\label{le42}
Let $u >0 $ be a  fixed real number. Then the  entire function \\
$f_u(s) := \txi_{\QQ}(u, s) : = \frac{1}{2} \frac{s(s-u)}{u} \zq (u,s)$ 
has all its zeros in the vertical strip
\begin{equation}\label{S403A}
|\Re (s) - \frac{u}{2}|< \frac{u}{2} +  8.
\end{equation}
\end{lemma}

\paragraph{Proof.}
We have
\begin{equation}\label{S405}
f_u (s) = \frac{s(s-u)}{4u} \pi^{- s/2} 
\Gamma \left( \frac{s}{2} \right) D_u \left( \frac{s}{2} \right)
\end{equation}
where
\begin{equation}\label{S406}
D_u (s) = \sum_{m=1}^\infty c_m (u) m^{-s} ~
\end{equation}
converges absolutely for $\Re(s) > \frac{u}{2} + 1$
by Theorem~\ref{th33}(i), and meromorphically continues to $\CC.$
All zeros of $f_u (s)$ with $\Re (s) > 0$ must come 
from those of the Dirichlet series $D_u ( \frac{s}{2} ) =0$.

The Dirichlet series $D_u (s)$ has no zeros in any half plane 
${\rm Re} (s) > \sigma$ for any $\sigma$ with 
$$
|c_1 (u) | > \sum_{m=2}^\infty |c_m (u) | m^{-\sigma} \,.
$$
Since $c_1 (u) = 2u$,  for $u > 0$ we may rewrite this as
\begin{equation}\label{S407}
1 > \frac{1}{2} \sum_{m=2}^\infty \frac{|c_m (u) |}{u} m^{-\sigma} 
\end{equation}
Now Theorem \ref{th33}(ii) gives
$$\frac{1}{2}\sum_{m=2}^\infty \frac{|c_m (u) |}{u} m^{-\sigma} \le 
3 \sum_{m=2}^\infty m^{- \sigma +\frac{u}{2}+1} ~.
$$
Choosing $\sigma = \frac{u}{2} + 4$ yields
$$\frac{1}{2}\sum_{m=2}^\infty \frac{|c_m (u) |}{u} m^{-\sigma} \le 
3 \sum_{m=2}^\infty m^{-3} = 3(\zeta(3) - 1) < 1,
$$
as required. Thus $D_u(\frac{s}{2})$ 
has no zeros in $\Re(s) > u + 8$, hence $f_u(s)$ also has no zeros there.
Finally the  functional equation $f_u(s) = f_u(u-s)$ 
implies that $f_u(s)$ has no zeros in the region
$\Re(s) < u - (u+8) = -8.$ ~~~$\bsq$

\noindent\paragraph{Remark.} The width of the 
strip  of Lemma~\ref{le42} is
qualitatively correct, in that it must grow like $u + O(1)$ for large $u$
and it must be of positive width, at least $4$, as $u \to 0$
to accomodate the zeros of $\xi_{\QQ}(0, s)$ given in Theorem~\ref{nth51}.

%
%
%

\subsection{Counting Zeros to Height $T$}

We establish the
following estimate for the number of zeros $N_u(T)$ within distance $T$ of the
real axis of $f_u(s)$, which generalizes a similar estimate for the
Riemann zeta function.

\begin{theorem}\label{SPth51}
There is an absolute constant $C_0$ such that, for all real $u > 0$ 
and $T \ge 0$,
\begin{equation}\label{SP511}
\frac{1}{2} N_u (T) = \frac{T}{2 \pi} \log \frac{T}{2 \pi} - \frac{T}{2 \pi} +
\frac{7}{8} + S_u (T),
\end{equation}
in which $S_u(T)$ satisfies
\begin{equation}\label{SP512}
|S_u (T) | \le C_0 (u + 1) \log (T+ u + 2).
\end{equation}
\end{theorem}

The proof of this result generalizes the proof for $\zeta(s)$
in Davenport \cite[Sect. 16]{Da80}, with some extra work to
control the dependence in $u$ in all estimates. 
We prove several preliminary lemmas.

We use the argument principle, and let $\Delta_L \arg (g(s))$ 
denote the change in argument $\theta$ in a function 
$g(s) = R(s) e^{i \theta(s)}$ along a contour $L$ on which 
$g(s)$ never vanishes.
For positive real $u$,
the zeros of $f_u(s)$ are those of 
the analytic continuation of the Dirichlet
series $D_u(\frac{s}{2})$ given in \eqn{S406}, possibly excluding
zeros of $D_u(\frac{s}{2})$  at negative even integers. 

We consider the rectangular contour $R$, oriented counterclockwise,
with corners at $\sigma_0 \pm iT$ and $u - (\sigma_0 \pm i T)$,
where
$$ \sigma_0 = u + 10, $$
and $T$ is chosen to avoid any zeros of $f_u(s).$ 
(See Figure 5.1) We will mainly
use the quarter-contour $L$ 
consisting of a vertical line $V$ from 
$z_0 = \sigma_0$ to $z_1 = \sigma_0 + i T$, 
followed by a horizontal line $H$ from $z_1$ to 
$z_2 = \frac{u}{2} + iT$.

\begin{lemma}\label{SPle52}
For real $u > 0$, and $T \ge 2$,
\begin{equation}\label{SP515}
\frac{1}{2} N_u (T) = \frac{T}{2 \pi} \log \frac{T}{2 \pi} - 
\frac{T}{2 \pi} + \frac{7}{8} + S_u (T)
\end{equation}
with
\begin{equation}\label{SP516}
S_u (T) = \Delta_L \arg D_u \left( \frac{s}{2} \right) + 
O \left( \frac{u}{T} \right) \,.
\end{equation}
\end{lemma}

\noindent\paragraph{Proof.}
We have
$$2 \pi N_u (T) = \Delta_R \arg ( f_u (s))$$
on the rectangular contour $R$ , oriented counterclockwise, 
which has its corners at $\sigma_0 \pm iT$ and 
$\frac{u}{2} \pm \sigma_0 \pm iT$,
as given above.

The functional equation $f_u (s) = f_u (u-s )$ 
and the symmetry $f_u (r) = \overline{f_u (\bar{s} )}$ imply that
\begin{equation}\label{SP517}
\frac{\pi}{2} N_u (T) = \arg_L (f_u(s)) \,,
\end{equation}
on the quarter-contour $L$, with each other quarter-contour 
of $R$ contributing the same amount.
Now
\begin{equation}\label{SP518}
u f_u (s) = \frac{1}{2}
(s-u) \pi^{- \frac{s}{2}} \Gamma \left( \frac{s}{2} +1 \right)
D_u \left( \frac{s}{2} \right) \,,
\end{equation}
so we obtain
\begin{equation}\label{SP519}
\frac{\pi}{2} N_u (T) =
\Delta_L \arg \left( \pi^{- \frac{s}{2}} \right) + 
\Delta_L \arg \left( \frac{s}{2} \Gamma (\frac{s}{2}) \right) +
\Delta_L \arg \left( s-u \right) + 
\Delta_L \arg \left( D_u(\frac{s}{2}) \right) \,.
\end{equation}
The first three terms on the right contribute
\begin{eqnarray}\label{SP520}
\Delta_L \arg \left( \pi^{- \frac{s}{2}}\right)  & = & 
\Delta_L \left( - \frac{1}{2} t \log \pi \right) 
= - \frac{T}{2}  \log \pi \,, \nonumber \\
\Delta_L \arg \left( s-u \right) & = & \arg \left( iT - \frac{u}{2} \right) =
\frac{\pi}{2}  + O \left( \frac{u}{T} \right) \,, \\
\Delta_L \arg \Gamma \left( \frac{s}{2} +1 \right) & = & \frac{T}{2}  
\log\left( \frac{T}{2}  \right) - \frac{T}{2}  
+ \frac{3\pi}{8}  + O \left( \frac{1}{T} \right) \nonumber\,
\end{eqnarray}
where Stirling's formula is used for the last estimate.
This yields \eqn{SP516}.~~~$\bsq$

Our object will be to  estimate $\Delta_L \arg D_u \left( \frac{s}{2} \right)$ 
using the formula
$$\Delta_L \arg \left( D_u ( \frac{s}{2}) \right) = 
- \int_L \Im \left( \frac{1}{2} \frac{D'_u ( \frac{s}{2})}
{D_u ( \frac{s}{2})} \right) ds \,,
$$
starting from the endpoint $s = \sigma_0$ of $L$ ,
where the next lemma shows $D_u \left( \frac{\sigma_0}{2} \right)$ is real
and positive
and 
$- \frac{1}{2} \frac{D'_u \left( \frac{\sigma_0}{2} \right)}
{D_u \left( \frac{\sigma_0}{2} \right)}$ 
is real and positive. In the integral we  analytically continue
$\frac{D'_u \left( \frac{s}{2} \right)}{D_u \left( \frac{s}{2} \right)}$
along $L$, and we choose $T$ so that the contour $L$ encounters no zero of
$D_u \left( \frac{s}{u} \right)$.
We need information on the zeros $\sZ_u$ of $D_u (s)$ 
obtained from the Hadamard product.

\begin{lemma}\label{le53}
Let $u > 0$ be real.
\begin{itemize}
\item[(i)]
For $s = \sigma + it \not\in \sZ_u \cup \{ -2,  -4,  - 6, \ldots \}$ 
there holds
\begin{equation}\label{SP521}
\sum_{\rho \in \sZ_u} \frac{\sigma - \Re(\rho)}{|s- \rho |^2} =
  \frac{1}{2} \Re \left( \frac{D'_u ( \frac{s}{2} )}
{D_u ( \frac{s}{2})} \right) +
 \frac{1}{2}\Re \left( \frac{\Gamma ' ( \frac{s}{2} +1 )}
{\Gamma ( \frac{s}{2} +1 )} \right) +
\frac{\sigma -u}{|s-u|^2}
- \frac{1}{2} \log \pi \,.
\end{equation}
\item[(ii)]
For $\Re (s) > u + 8$,
\begin{equation}\label{SP522}
\Re \left( D_u \left( \frac{s}{2} \right) \right) > 0 \,,
\end{equation}
and there is an absolute constant $A_0$ independent of $u$ such that
\begin{equation}\label{SP523}
\left| \frac{D'_u \left( \frac{s}{2} \right)}
{D_u \left( \frac{s}{2} \right)} \right| \le A_0.
\end{equation}
\end{itemize}
\end{lemma}

\paragraph{Proof.}(i). We use the Hadamard product expansion
\begin{equation}\label{SP524}
u f_u (s) = \frac{1}{2}\frac{s(u-s)}{u} Z(u,s) = e^{A(u) + B(u) s}
\prod_{\rho \in \sZ_u} \left( 1 - \frac{s}{\rho} \right) e^{\frac{s}{\rho}}
\end{equation}
where $\sZ_u$ is the set of zeros of $f_u (s)$ counted with multiplicity.
 This formula  is valid 
because $u f_u (s)$ is an entire
function of order at most one by the growth estimate of Theorem \ref{Nth41}.
Note that $0 \not\in \sZ_u (s)$ because for fixed $u$ the function
$Z(u,s)$ has a 
simple pole with residue 1 at $s=0.$ 
The derivation of \cite[pp. 82--84]{Da80} yields the formula 
\begin{equation}\label{SP528A}
B(u) = \Re (B(u)) = - \sum_{\rho \in \sZ_u} \!^{\!^\prime}~ \frac{1}{\rho} = 
- \sum_{\rho \in \sZ_u} \frac{\Re(\rho)}{|\rho |^2} \,,
\end{equation}
where the prime in the first sum indicates that complex conjugate
zeros $\rho$ and $\bar{\rho}$ are to be summed in pairs, and the
last sum converges absolutely.

We set equal the logarithmic derivatives of \eqn{SP518} and \eqn{SP524}, 
to obtain
\begin{equation}\label{SP524a}
\frac{1}{2} \frac{D'_u \left( \frac{s}{2} \right)}
{D_u \left( \frac{s}{2} \right)} + 
\frac{1}{2} \frac{\Gamma ' \left( \frac{s}{2} +1 \right)}
{\Gamma \left( \frac{s}{2} +1 \right)} +
\frac{1}{s-u} - \frac{1}{2} \log \pi 
 =  B(u) + \sum_{\rho \in \sZ_u} \left( \frac{1}{s- \rho} + 
\frac{1}{\rho} \right) \,.
\end{equation}
This yields
\begin{equation}\label{SP525}
- \frac{1}{2} \frac{D'_u \left( \frac{s}{2} \right)}
{D_u \left( \frac{s}{2} \right)} = 
\frac{1}{s-u} -  B(u) - \frac{1}{2} \log \pi 
+ \frac{1}{2} \frac{\Gamma' \left (\frac{s}{2} +1 \right)}
{\Gamma \left( \frac{s}{2} +1 \right)}
- \sum_{\rho \in \sZ_u}
\left( \frac{1}{s-\rho} + \frac{1}{\rho} \right). 
\end{equation}
Taking real parts yields,
\begin{eqnarray}
\label{SP526}
- \frac{1}{2} \Re
\left( \frac{D'_u \left( \frac{s}{2} \right)}
{D_u \left( \frac{s}{2} \right)} \right) & = &
\frac{\sigma- u}{|s-u|^2} - \Re (B(u)) - \frac{1}{2} \log \pi \nonumber \\
&& \quad + \, \frac{1}{2} \Re \left(
\frac{\Gamma ' \left( \frac{s}{2} +1 \right)}
{\Gamma' \left( \frac{s}{2} +1 \right)} \right) 
- \sum_{\rho \in \sZ_u} \left( \frac{\sigma- \Re(\rho)}{|s - \rho |^2} +
\frac{\Re(\rho)}{|\rho|^2} \right) \,.
\end{eqnarray}
Now \eqn{SP526} holds in the entire plane by analytic continuation,
since the functions are single-valued.
Applying the formula for $B(u)$ simplifies 
\eqn{SP526} to \eqn{SP521}, proving (i).

(ii) Suppose $\Re (s) > u + 8$.
Then the formula \eqn{SP518} for $D_u (s)$ and Theorem \ref{th33}(ii) give
\begin{eqnarray}\label{SP527A}
\left| D_u \left( \frac{s}{2} \right) - 2u \right| & \le &
\sum_{m=2}^\infty |c_m (u) | m^{- \sigma/2} \nonumber \\
& \le & 6u \sum_{m=1}^\infty m^{\frac{u +1 -\sigma}{2}} \le 6u(\zeta(3) - 1)
 \nonumber \\
& \le & \frac{3}{2} u \,.
\end{eqnarray}
This implies \eqn{SP522}.
Next, applying Theorem \ref{th33}(ii) again,
\begin{eqnarray*}
\left| \frac{D'_u \left( \frac{s}{2} \right)}
{D_u \left( \frac{s}{2} \right)} \right| & \le & 
\frac{1}{\left(2 - \frac{3}{2} \right) u} 
\left| D'_u \left ( \frac{s}{2} \right) \right| \\
& \le & \frac{2}{u} \sum_{m=2}^\infty 
|c_m (u) | \frac{\log m}{2} m^{- \sigma/2} \\
& \le & 
24 \sum_{m=2}^\infty (\log m) m^{\frac{u + 1 - \sigma}{2}} \\
& \le & 
24 \sum_{m=2}^\infty (\log m) m^{- 3} = O(1) \,,
\end{eqnarray*}
which gives \eqn{SP523}.~~~$\bsq$

\begin{lemma}\label{le54}
(i) There is an absolute constant $A_1$, such that for all $u > 0$, 
and  all real $T$,
\begin{equation}\label{SP527}
\sum_{\rho = \beta + i \gamma \in \sZ_u} 
\frac{u + 2}{4(u+ 9 )^2 + (T- \gamma)^2}
\le 
A_1 \log (|T| + u + 2 ) \,.
\end{equation}

(ii) There is an absolute constant $A_2$ such that for all $u > 0$ and all $T$,
the number of zeros $\rho = \beta + i \gamma \in \sZ_u$  with
\begin{equation}\label{SP528}
| \gamma - T | < u + 9,
\end{equation}
counting multiplicity, is at most
\begin{equation}\label{SP529}
A_2 (u+ 9 ) \log ( |T| + u + 2 ).
\end{equation}
\end{lemma}

\paragraph{Proof.}
(i). Choose $s = \sigma_0 + iT$ with
\begin{equation}\label{SP530}
\sigma_0 = 2u + 10 \,,
\end{equation}
so 
$\frac{\sigma_0 - u}{|s - u|^2} \le  \frac{u + 10}{|\sigma_0 - u|^2} 
\le \frac{1}{10}.$
Now apply \eqn{SP521} and the bound \eqn{SP523} to obtain
\begin{equation}\label{SP531}
\sum_{\rho \in \sZ_u} \frac{\sigma_0 - \Re{\rho}}{| s- \rho |^2} = 
 \frac{1}{2} \Re \left( \frac{\Gamma ' \left( \frac{s}{2} +1 \right)}
{\Gamma \left( \frac{s}{2} +1 \right)} \right)  + O (1)\,.
\end{equation}
We recall the formula
\begin{equation}\label{SP532}
\frac{\Gamma' (s)}{\Gamma (s)} = \log s + O \left( \frac{1}{|s|} \right)
\end{equation}
valid for $- \pi + \delta < | \arg s | < \pi - \delta$ for any fixed $\delta$.
(We choose $\delta = \frac{\pi}{2}$.)
Now \eqn{SP530} gives
\begin{eqnarray}
\label{SP533}
\Re \left( \frac{\Gamma ' \left( \frac{s}{2} +1 \right )}
{\Gamma \left( \frac{s}{2} + 1 \right)} \right) 
& = & 
\Re \left( \log \frac{s}{2} +1 \right) + O(1) \nonumber \\
& = & 
\log \left| \frac{s}{2} +1 \right| + O(1) \nonumber \\
& = & 
\log ( | \sigma_0 | + |T| ) \,.
\end{eqnarray}
Thus \eqn{SP531} becomes
\begin{equation}\label{SP534}
\sum_{\rho \in \sZ_u} 
\frac{\sigma_0 - \beta}
{(\sigma_0 - \beta )^2 + (T- \gamma)^2} =
O( \log (|T| + 2u + 8)) ~.
\end{equation}
The bound of Lemma \ref{le42} gives
$$u + 2 \le \sigma_0 - \Re(\rho) \le 2(u + 9) ~.$$
Thus
$$\frac{u+ 2}{4(u+ 9 )^2 + 
(T- \gamma )^2} \le \frac{\sigma_0 - \beta}
{(\sigma_0 - \beta)^2 + (T- \gamma)^2} ~,
$$
and \eqn{SP527} follows.

(ii). Let $S_T$ denote the set of zeros in $\sZ_u$ 
satisfying $| \gamma - T | < u + 9.$
Then, since $u \geq 0$,
$$\sum_{\rho \in S_T}
\frac{u+ 2}{4(u+ 9 )^2 + | \gamma - T |^2} 
\ge 
\frac{2}{45}
\left(\frac{|S_T|}{u+ 9} \right)~.
$$
Combining this with \eqn{SP527} implies \eqn{SP529}
 with $A_2 = \frac{45}{2} A_1$.~~~$\bsq$

\begin{lemma}\label{le55}
There is an absolute constant $A_3$, such that for $u > 0$ and 
$s = \sigma + iT$ with $s \not\in \sZ_u$ in the region
\begin{equation}\label{SP535}
\left| \sigma - \frac{u}{2} \right| \le \frac{3u}{2} + 10 ~,
\end{equation}
there holds, for all  $|T| \ge |\sigma | + 1,$  
\begin{equation}\label{SP537}
\left| \frac{1}{2} \frac{D'_u \left( \frac{s}{2} \right)}
{D_u \left( \frac{s}{2} \right)} - 
\sum_{\substack{\rho \in \sZ_u\\| \gamma - T | \le u+9}} \frac{1}{s-\rho} \right|
\le A_3 \log ( |T| + u + 2 ).
\end{equation}
\end{lemma}

\paragraph{Proof.}
Set $\sigma_0 = 2u + 10$, and
$s_0 = \sigma_0 + iT$.
Then,  differencing \eqn{SP525} at $s$ and $s_0$, we obtain
\begin{eqnarray*}
\frac{1}{2}
\frac{D'_u \left( \frac{s}{2} \right)}{D_u \left( \frac{s}{2} \right)} & - &
\frac{1}{2}
\frac{D'_u \left( \frac{s_0}{2} \right)}{D'_u \left( \frac{s_0}{2} \right)} =
\left( - \frac{1}{s-u} + \frac{1}{s_0 -u} \right) \\
&& \quad + \,
\sum_{\rho \in \sZ_u}
\left\{\left( \frac{1}{s-\rho} +
\frac{1}{\rho} \right) -
\left( \frac{1}{s_0 - \rho} + \frac{1}{\rho} \right)\right\} \\
& - & \frac{1}{2}
\left( 
\frac{\Gamma ' \left( \frac{s}{2} +1 \right)}
{\Gamma' \left( \frac{s}{2} +1 \right)} -
\frac{\Gamma' \left( \frac{s_0}{2} +1 \right)}
{\Gamma' \left( \frac{s_0}{2} +1 \right)} \right) \,.
\end{eqnarray*}
The bound $|T| \ge |\sigma | + 1$  
ensures that $s = \sigma + iT$ has 
$- \pi + \delta < \arg (s) < \pi + \delta$ 
for $\delta = \frac{\pi}{4}$, hence the bounds
\eqn{SP533} applies to give
\begin{equation}\label{SP538}
\frac{1}{2}
\frac{D'_u \left( \frac{s}{2} \right)}{D_u \left( \frac{s}{2} \right)}  =
\sum_{\rho \in \sZ_u} \left( \frac{1}{s-\rho} - 
\frac{1}{s_0 - \rho} \right) + O(1) +  O ( \log (|T| + u + 2) )\,.
\end{equation}
Now
\begin{eqnarray}\label{SP539}
\left| \frac{1}{s-\rho} - \frac{1}{s_0 - \rho} \right|
=
\frac{|s-s_0|}{|s- \rho | | s_0 - \rho |} & \le &
\frac{3u + 20}{|s- \rho | | s_0 - \rho |} \nonumber \\
& \le & \frac{4(u + 5)}{| \gamma - T|^2} \,.
\end{eqnarray}
For those zeros with $|\gamma - T | > u + 9$, 
Lemma \ref{le54}(i) gives the bound
\begin{eqnarray*}
\sum_{\substack{\rho \in \sZ_u \\ |\gamma -T > u + 9}}
\frac{4(u + 5)}{(\gamma -T)^2} 
& \le & 
60\sum_{\rho \in \sZ_u} \frac{u + 2}{4(u + 9)^2 + | \gamma -T|^2} \nonumber \\
& \le & 
60A_1 ( \log ( |T| + u + 2 )) .
\end{eqnarray*}
If $| \gamma - T | \le u + 9$, then by Lemma \ref{le54}(ii) 
the number of such zeros is $O((u + 1) \log ( |T| + u + 2 ))$ 
and for each one
$$\left| \frac{1}{s_0 - \rho} \right| =
O \left( \frac{1}{|\sigma_0 - \beta|} \right) = 
O \left( \frac{1}{u + 9} \right) \,.
$$
So their total contribution is $O( \log |T| + u^\ast ))$ in \eqn{SP538}.
Substituting these bounds in \eqn{SP538} yields
$$\left| \frac{1}{2} \frac{D'_u \left( \frac{s}{2} \right)}
{D_u \left( \frac{s}{2} \right)} -
\sum_{\substack{\rho \in \sZ_u \\ | \rho - T | < u + 9}}
\frac{1}{s-\rho} \right| =
O( \log ( |T| + u + 2 )) \,,
$$
as required.~~~$\bsq$

\paragraph{Proof of Theorem \ref{SPth51}.}
Set $\sigma_0 = 2u + 10$.
Recall that the quarter- contour $L$ consists of the
 vertical segment $V$ from $\sigma_0$ to $\sigma_0 + iT$, and 
the horizontal segment $H$ from $\sigma_0 + iT$ to $\frac{u}{2} + iT$,
so that 
\begin{equation}\label{SP540}
\pi S_u ( T) = \Delta_V \arg D_u \left( \frac{s}{2} \right) + 
\Delta_H \arg D_u \left( \frac{s}{2}\right) + O\left( \frac{u}{T} \right)
\end{equation}
by  \eqn{SP520}.
We have
$$\left| \Delta_V \arg D_u \left( \frac{s}{2} \right) \right|
\le \pi,$$
using Lemma \ref{le53}(ii). Now
$$\Delta_H \arg D_u \left( \frac{s}{2} \right) = 
- \int_{\frac{1}{2} + iT}^{\sigma_0 + iT} \Im 
\left(\frac{1}{2} \frac{D'_u \left( \frac{s}{2} \right)}
{D_u \left( \frac{s}{2} \right)} \right) ds \,.
$$
and to estimate this we
apply Lemma \ref{le55}. For each zero $\rho$ we have
$$\int_{\frac{1}{2} + iT}^{\sigma_0 + iT}
\Im \left( \frac{1}{s- \rho} \right) ds = \Delta_H \arg (s- \rho ),
$$
which contributes  at most $\pi$.
Now Lemma \ref{le54}(ii) gives 
that there are at most
$O((u + 1) \log (|T| + u + 2 ))$ such zeros in the sum \eqn{SP537}, 
so their total contribution to the argument is  at most
$O((u + 1) \log ( |T| + u + 2 ))$.
The error term in \eqn{SP537} integrated over 
$H$ contributes at most a further
$O( (u + 1)\log (|T| + u + 2 ))$ to the argument, since
$H$ is a  path of length at most $\frac{3u}{2} + 10$.
Combining these estimates in \eqn{SP540} yields the bound \eqn{SP512}, 
completing the proof.~~~$\bsq$

\noindent\paragraph{Remark.} The zero-counting estimate of Theorem~\ref{SPth51}
for $u > 0$ is easily checked to remains valid for $u =0$ by virtue of 
the explicit formula for $f_0(s)= \xi_{\QQ}(0, s)$ in Theorem~\ref{nth51}.
Note that the extra zeros provided by the terms 
$(1 - 2^{1 + \frac{s}{2}})(1 - 2^{1 - \frac{s}{2}})$ are needed to make
the main term \eqn{SP511} valid.

%
%
%

\subsection{Movement of Zeros}

It is interesting to examine the behavior of the zero set
of $\xi_{\QQ}(u, s)$ for nonnegative real $u$, as $u$ is varied. 
As above, consider  variation $1 \le u \le 2$, or, more generally, over a 
fixed bounded range of $u$. For that range of $u$, Theorem~\ref{SPth51}
asserts that the  general density of zeros to height $T$
remains almost constant, with a variation of $O(\log (T+2)).$
Since the number of zeros in a unit interval at this height is
of the same order, it suggests that
every zero can move vertically
a distance of at most $O(1)$, while  the horizontal movement is 
certainly  restricted
to distance $O(1)$ by Lemma~\ref{le42}. Therefore it
 would appear that varying $1 \le u \le U$, there is a constant
$C_U$ depending on $U$ such that every zero moves at most 
the bounded amount $C_U$, {\em independent of the height $T$ of this zero} 
in the critical strip.  We cannot assert this rigorously,
however, because we have not ruled out the possibility of zeros going
off the line in pairs and then hop-scotching around other zeros
remaining on the line.

\begin{table}[p]
$$
\begin{array}{|r|r|rr|}
\hline
&\multicolumn{1}{c|}{\zeta(s)}&\multicolumn{2}{c|}{\zeta_{\QQ(i)}(s/2)}\\
\cline{2-4}
 1  ~&~  14.13  ~&~  12.04  &~\\
 2  ~&~  21.02  ~&~  20.48  &~\\
 3  ~&~  25.01  ~&~  25.96  &~\\
 4  ~&~  30.42  ~&~  28.26  &*~\\
 5  ~&~  32.94  ~&~  32.68  &~\\
 6  ~&~  37.58  ~&~  36.58  &~\\
 7  ~&~  40.91  ~&~  42.04  &*~\\
 8  ~&~  43.32  ~&~  42.90  &~\\
 9  ~&~  48.00  ~&~  46.54  &~\\
10  ~&~  49.77  ~&~  50.02  &*~\\
11  ~&~  52.77  ~&~  51.44  &~\\
12  ~&~  56.44  ~&~  56.78  &~\\
13  ~&~  59.34  ~&~  59.30  &~\\
14  ~&~  60.83  ~&~  60.85  &*~\\
15  ~&~  65.11  ~&~  65.18  &~\\
16  ~&~  67.07  ~&~  65.87  &*~\\
17  ~&~  69.54  ~&~  68.38  &~\\
18  ~&~  72.06  ~&~  72.28  &~\\
19  ~&~  75.70  ~&~  75.16  &*~\\
20  ~&~  77.14  ~&~  77.02  &~\\
21  ~&~  79.33  ~&~  80.64  &~\\
22  ~&~  82.91  ~&~  81.82  &*~\\
23  ~&~  84.73  ~&~  83.60  &~\\
24  ~&~  87.42  ~&~  86.64  &*~\\
25  ~&~  88.81  ~&~  89.22  &~\\
\hline
\end{array}
$$
\caption{Imaginary part of zeros of $\zeta(s)$ and $\zeta_{\QQ(i)}(s/2)$.
 ($*=\text{zero of }\zeta(s/2)$.)}
~\label{zeros}
\end{table}

Theorem ~\ref{SPth51} for $u=1$ and $u=2$ shows that the zero-counting
functions of $\zeta(s)$ and $\zeta_{\QQ(i)}(\frac{s}{2})$ have 
extemely similar asymptotics. 
In Table~\ref{zeros} we compare the first $25$ such
zeros. Since $\zeta_{\QQ(i)}(s)=\zeta(s) L(s, \chi_{-4})$, we
have indicated the zeros of $\zeta_{\QQ(i)}(\frac{s}{2})$ associated to
$\zeta(\frac{s}{2})$ with an asterisk, and the remainder
come from $L(\frac{s}{2}, \chi_{-4})$. Note that the appearance
of every zeta zero on both sides of this table shows that
zeta zeros have a kind of self-similar structure by powers of two.
However this is only approximately true, because
there is  no precise correspondence of zeta zeros due
to the phenomonon of zeros coalescing as $u$ varies, as indicated
in Lemma~\ref{Nle71} in \S8. 

%
%
%

\section{Location of Zeros: $\cw$ Negative Real}
\hsp
In this section we  suppose $\cw = -u$  is negative real 
($u > 0$) and
fixed. We will not find zeros, but will instead
 specify places in the $s$-plane where the the function
$\xi_{\QQ}(u, s)$ has no zeros. An interesting extra structure 
underlying certain properties of  the two-variable zeta function 
$\zq(\cw, s)$ for negative real $\cw$ 
is a holomorphic convolution semigroup of 
complex-valued measures
described in \S7.2. On the critical line $\Re(s) = - \frac{u}{2},$
we will show these are real-valued positive measures, normalizable
to be probability measures, in \S7.1.

%
%
%

\subsection{Positivity on Critical Line}

We  will establish the following result, concerning
the absence of zeros on the ``critical line'' $\Re(s) = -\frac{u}{2}$
($u >0$), which will be deduced from Theorem~\ref{nth62} below.

\begin{theorem}\label{nth61}
$($Positivity Property$)$
For real  $\cw = -u$ with  $u > 0$ and all 
real $t$,
\begin{equation}\label{n61}
\zq \left( -u, -\frac{u}{2} + it \right) > 0 ~.
\end{equation}
\end{theorem}

Recall that the functional equation on the ``critical line'' 
$s = -\frac{u}{2} + it$ implies that
$$
\zq \left( -u, -\frac{u}{2} + it \right) = 
\zq \left( -u,-\frac{u}{2} -it \right) = 
\overline{\zq \left( -u, -\frac{u}{2} + it \right)} ~,
$$
hence $\zq(-u, -\frac{u}{2} + it )$ is real.
For real $-u < 0$, the integral representation
\begin{equation}~\label{n61a}
\zq \left( -u, -\frac{u}{2} + it \right) = \int_0^\infty \theta (x^2)^u
x^{-\frac{u}{2} + it } \frac{dx}{x}
\end{equation}
converges absolutely for all $t \in \RR$, and it gives
\begin{equation}\label{n62}
\zq \left( -u, -\frac{u}{2} \right) > 0 ~,
\end{equation}
since the integrand is positive.
Thus the assertion that
$\zq (-u, -\frac{u}{2} + it ) \ne 0$ for all real $t$ 
is equivalent to the positivity condition \eqn{n61}.
By changing  variables in the integral \eqn{n61a}, with $x = e^{-r}$, we obtain
\begin{eqnarray}\label{n63}
\frac{1}{2\pi}\zq\left( -u, -\frac{u}{2} + it \right) & = &
\frac{1}{2\pi}\int_{-\infty}^\infty 
\theta (e^{-2r})^{-u} e^{\frac{ru}{2}} e^{-irt} dr
\nonumber \\
& = & \frac{1}{2\pi}\int_{-\infty}^\infty
\left( \frac{e^{r/2}}{\theta (e^{-2r})} \right)^{u} e^{-irt} dr. 
\end{eqnarray}
This implies that the
function
\begin{equation}\label{n66}
P_u(x) := \frac{1}{2\pi} \theta(1)^u Z_{\QQ}(-u, -\frac{u}{2} +ix)
\end{equation}
has inverse Fourier transform
\begin{equation}\label{n64b}
\check{P}_u(r) := \int_{-\infty}^\infty P_u(x) e^{ixr} dx
= \theta(1)^u \left( \frac{e^{r/2}}{\theta (e^{-2r})} \right)^{u}.
\end{equation}  
The function $P_u(x)$ is an even function and has
\begin{equation}~\label{n64c}
\int_{-\infty}^\infty P_u(x)dx = \check{P}_u(0)=1.
\end{equation}
This equation shows that $P_u(x)dx$ is a signed measure of
mass one. 
Theorem \ref{nth61} asserts that 
for $P_{u} (x) > 0$ for all $x$ holds for  all $u >0$,
which would imply that $P(x)dx$
is a probability measure.
In any case $\check{P}_u(r)$ is the characteristic function of
the (signed) measure $P_{u} (x) dx$, and \eqn{n64b} shows that  
$\check{P}_u(r) = f(r)^u$ where
\begin{equation}\label{n65}
f(r) :=
 \theta (1) \frac{e^{r/2}}{\theta (e^{-2r})} =
 \frac{\theta (1)}{\sqrt{\theta (e^{2r} ) 
\theta (e^{-2r} )}},
\end{equation}
where the last expression is derived using the functional equation for
the theta function.
The assertion that $f(r)^{u}$ is a characteristic function for
all real $u > 0$ is equivalent to the assertion that each
$P_u (x) dx$ is an {\em infinitely divisible} probability measure;
the collection $\{P_u (x) dx ~: u > 0\}$ then form a 
semigroup under convolution.

We note that the normalizing factor
\begin{equation}\label{n68}
\theta (1) = 1+ 2 \sum_{n=1}^\infty e^{- \pi n^2} =
\frac{\pi^{\frac{1}{4}}}{\Gamma \left( \frac{3}{4} \right)}
\approx 1.08643
\end{equation}
appearing in \eqn{n66} is the invariant $\eta (\QQ )$ 
introduced by van~der Geer and Schoof \cite[p. 16]{vdGS99}.
They define the {\em genus of $\QQ$}
to be the ``dimension'' $h^0 (\sK_{\QQ} )$ of the canonical divisor
$ \sK_{\QQ}  = (1)$, which is
\begin{equation}\label{n69}
h_0 ( \sK_{\QQ} ) = \log
( \eta ( \QQ ) \sqrt{\Delta_{\QQ}}) = \log \theta (1) ~,
\end{equation}
see the appendix.

A necessary and sufficient condition for a function to be 
a characteristic function of an infinitely divisible probability
measure was developed by Khintchine, Levy and Kolmogrov.
We follow the treatment in Feller \cite[p. 558--563]{Fe71}.
A measure $M\{dy\}$ on $\RR$ is called {\em canonical} if it is nonnegative, 
assigns finite masses to finite intervals and if both the integrals
\begin{equation}\label{n610}
M^+ (x) = \int_x^\infty
\frac{1}{y^2} M \{dy\} , ~~
M^- (x) = \int_{- \infty}^{-x} \frac{1}{y^2} 
M \{dy\} 
\end{equation}
converge for some (and therefore all) $x > 0$.

\begin{prop}\label{npr61}
A complex-valued function $f(r)$ on $\RR$ is the characteristic function 
of an infinitely divisible probability measure if and only if 
$f(r) = \exp(\psi (r))$ with $\psi(r)$ having the form
\begin{equation}\label{n612}
\psi (r) = ibr + \int_{-\infty}^\infty
\frac{e^{irx} -1-ir \sin x}{x^2}
M \{dx\}
\end{equation}
for some canonical measure $M$ and real constant $b$.
The canonical measure $M$ and constant $b$ are unique.
\end{prop}

\paragraph{Proof.}
This is shown in  Feller \cite[pp. 558--563]{Fe71}.~~~$\bsq$

The representation \eqn{n612}
implies that $f(0)=1$, $f(-r)=\overline{f(r)}$ for all $r$ 
and that $\log f(r) $ is well-defined, with its imaginary part
determined by continuity starting from $\log f(0) = 0.$

We call the measure $M \{ dx \}$ in Proposition~\ref{npr61}, which may have 
infinite mass, the {\em Feller canonical measure} associated to $f(r)$.
A related canonical measure is the 
{\em Khintchine canonical measure} $K\{dx\}$,  given by
$$K\{dx\} = \frac{1}{1+x^2} M \{dx \}.$$ 
It can 
be an  arbitrary bounded nonnegative measure, see
Feller~\cite[pp. 564-5]{Fe71}. 
Biane, Pitman and Yor~\cite[p. 9]{BPY99} consider 
the {\em Levy-Khintchine
canonical measure} $\vu \{ dx \},$
which is
defined for infinitely divisible distributions supported on $[0, \infty)$,
and is  related to the corresponding Feller canonical measure by 
$$
\vu\{dx\} = \frac{1}{x^2} M \{dx \}.
$$
Note that the Feller canonical measure $M  \{ dx \}$ 
given in Theorem~\ref{th11}
has support on the whole real line, so has no associated
Levy-Khintchine canonical measure.

We will use a variant of this result which characterizes 
infinitely divisible distributions with finite second moment.

\begin{prop}\label{npr62}
A complex-valued function $f(r)$ on $\RR$ is the characteristic 
function of an infinitely divisible probability distribution 
having a finite second moment if
(and only if) $f(r)$ is a $C^2$-function, $f(0)=1$, $f(-r)=\overline{f(r)}$,
for all $r$, $f(r) \ne 0$, and 
\begin{equation}\label{n613}
g(r) := \frac {f''(r) f(r) - f'(r)^2}{f(r)^2} = \frac{d^2}{dr^2} ( \log f(r))
\end{equation}
has $g(0) < 0$ and $-g(r)$ is the 
characteristic function of a positive measure $M \{dx \}$ of 
finite mass,
\begin{equation}\label{n614}
-g(r) = \int_{-\infty}^\infty e^{irx} M \{dx\} ~.
\end{equation}
If so, then $M\{dx\}$ is the associated Feller canonical measure to
$f(r)$.
\end{prop}

\paragraph{Proof.}
The ``if'' part of this  result appears in Feller \cite[p. 559 bottom]{Fe71}.
We will not use the ``only if'' part of the result and omit its
proof.~~~$\bsq$

\begin{theorem}\label{nth62}
The function
\begin{equation}\label{n615}
f(r) = \theta (1)\frac{e^{r/2}}{\theta (e^{-2r} )} 
=\frac{\theta (1)}{\sqrt{\theta (e^{-2r}) \theta (e^{2r} )}} 
= \theta (1)\frac{e^{-r/2}}{\theta (e^{2r} )}~,
\end{equation}
on $\RR$ is the characteristic function of an infinitely
divisible probability measure with finite second moment.
Its associated Feller canonical measure
$M\{dx\}$ is equal to $M(x) dx$ with
\begin{equation}\label{n616}
M(x) = \frac{1}{8\pi}x^2
| 1 - 2^{1 + \frac{ix}{2}}|^2 
\left| \hat{\zeta} \left( \frac{ix}{2} \right) \right|^2.
\end{equation}
in which $\hat{\zeta} (s) := \pi^{- \frac{s}{2}} \Gamma \left(
\frac{s}{2} \right) \zeta (s)$. 
\end{theorem}

\noindent\paragraph{Proof.}
We apply the criterion of Proposition~\ref{npr62}. We start from the function
$$
g(r) = \frac{d^2}{dr^2} log (\theta(1) \frac{e^{-r/2}}{\theta(e^{2r})}) 
= - \frac{d^2}{dr^2}  \log \theta(e^{2r}), 
$$
and must show that $g(0) < 0$ and that 
$$
M(x) := - \frac{1}{2\pi} \int_{-\infty}^\infty e^{-irx} g(r) dr 
$$
is a nonnegative function having finite mass. We have
$g(0) \approx -1.8946 < 0$ and 
$$
M(x) =\frac{1}{2\pi} \int_{-\infty}^\infty e^{-irx}\left(\frac{d^2}{dr^2}  \log \theta(e^{2r})\right)dr= 
\frac{1}{2\pi} \int_{-\infty}^\infty e^{-irx} \gamma(0, r) dr,
$$
using \eqn{eq44a}. Now Theorem~\ref{th44} gives, on
taking $w=0$ and $y=0$, that 
for all $s \in  \CC$,
$$
\xi_{\QQ}(0, s)= \frac{1}{2} \int_{-\infty}^\infty e^{rs} \gamma(0, r) dr,
$$
and it follows that 
$M(x) = \frac{1}{\pi} \xi_{\QQ}(0, -ix).$
Now Theorem~\ref{nth51},
which uses the Jacobi triple product formula, 
gives
$$
\xi_{\QQ}(0, -ix) = \xi_{\QQ}(0, ix) =
\frac{x^2}{8} (1 - 2^{1 + \frac{ix}{2}})
(1 - 2^{1 + \frac{-ix}{2}})\hat{\zeta}(\frac{ix}{2}) 
\hat{\zeta}(\frac{-ix}{2}),
$$
Using
$\hat{\zeta}(\frac{-ix}{2})= \overline{ \hat{\zeta}(\frac{ix}{2}) }$
we obtain
$$
M(x) = \frac{1}{\pi}\xi_{\QQ}(0, ix) = \frac{1}{8\pi} x^2 | 1 - 2^{1 + \frac{ix}{2}}|^2
|\hat{\zeta}(\frac{ix}{2})|^2.
$$
This shows that $M(x)$ 
is nonnegative, and its strict positivity follows from
the well-known result that $\zeta(s)$ is
nonzero on the line $\Re(s)=1$, using 
$\hat{\zeta}(ix) = \hat\zeta(1 - ix)$.
The positive  measure
$M(x){dx}$ has finite mass since
$\xi_{\QQ}(0, ix)$ is a Schwartz function by Theorem~\ref{th43},
so the result follows by Proposition~\ref{npr62}.
(The mass of $M(x){dx}$ is explicitly determined in Theorem~\ref{thm79} below;
numerically it is about $1.8946.$)
~~$\bsq$

\noindent\paragraph{Proof of Theorem~\ref{nth61}.}
For real $u >0$ the function $f(r)^{u}$ with
$f(r)$ given by \eqn{n615} is by Theorem~\ref{nth62} the
characteristic function of a probability density of
a nonnegative measure. The measure $P_u (x) dx$ is
positive except on a discrete set because the density $P_u(x)dx$ is
an analytic function of $x$.  Then, using the infinite
divisibility property, we have
$$P_u(y) = P_{u/2} * P_{u/2} (y) = \int_{-\infty}^\infty P_{\frac{u}{2}}(x) P_{\frac{u}{2}}(y-x)dx 
>0,$$ 
giving positivity for all real $y$.
The probability density 
$P_u (x) dx$ is given 
by \eqn{n66}, and $\theta(1) > 0$,
so we conclude that 
$Z_{\QQ} (-u, -\frac{u}{2} + ix) > 0 $
for all real $x$.  
~~~$\bsq$

%
%
%

\subsection{Holomorphic Convolution Semigroup}

In \S7.1 we
showed that for $u > 0$ the family of probability measures  
$\{\rho_u := P_{u}(x) dx:~  u > 0\}$ having
the  density functions
$$
P_{u} (x) = \frac{1}{2\pi}\theta (1)^{u} 
\zq \left( -u, \frac{-u}{2} + ix\right) ~,
$$
form a semigroup under convolution, i.e.
$ \rho_{u_1} * \rho_{u_2} = \rho_{u_1 + u_2}.$
We can extend
this to a convolution semigroup of
complex-valued measures on the real line, 
indexed by two real parameters $(u, v)$ with $u > 0$ and $|v| < u$,
a region which forms an open cone in $\RR^2$ closed under addition.
 Given such $(u, v)$ we define
a complex-valued measure $\rho_{u, v}(x)dx$ on the real line by
\begin{equation}~\label{n718a}
\rho_{u, v}(x) =  \frac{1}{2\pi}\theta(1)^{u} 
\zq( - u, -\frac{u + v}{2} + ix) \qquad\mbox{for}\qquad -\infty < x < -\infty. 
\end{equation}
In terms of the function $\zq(\cw, s)$ these values occupy the 
real-codimension one cone 
\begin{equation}\label{n718b}
\sC^{-} = \{(\cw, s)~:~ \cw = u \in \RR ~~\mbox{and}~~ u  < \Re(s) < 0\}, 
\end{equation}
which
is contained in the absolute convergence region $\sC$ of the 
integral representation \eqn{S105}.

\begin{lemma}~\label{lem75}
For all $u>0$ and real $|v|<u$, $\rho_{u,v}(x)$ is a complex-valued
measure with
\[
\int_{-\infty}^\infty \rho_{u,v}(x) dx   = 1.
\]
The measures $\{ \rho_{u,v} : u > 0,~ |v|< u\}$ form a 
convolution semigroup. 
That is, for  
 $u_1, u_2>0$, $|v_1|< u_1, |v_2|< u_2$
and $x\in \RR$, we have:
\[
\int_{-\infty}^\infty \rho_{u_1, v_1}(y)\rho_{u_2, v_2}(x-y) dy
=
\rho_{u_1 + u_2, v_1+v_2}(x).
\]
In particular
\[
\int_{-\infty}^\infty \rho_{u,v}(x)^2 dx = \rho_{2u, 2v}(0).
\]
\end{lemma}

\noindent\paragraph{Proof.}
The first formula follows from the formula
\begin{equation}~\label{Nn638}
\zq(-u, -\frac{u+v}{2} + ix) =
 \check{\sF}\left( (\frac{e^r}{\theta(e^{-2r})})^{u} e^{vr}\right)
(-\frac{x}{\pi})
\end{equation}
which generalizes \eqn{n63}. Indeed the functions
$$
f_{u, v}(r) = \left(\frac{e^{\frac{r}{2}}}{\theta(e^{-2r})}\right)^{u} 
e^{\frac{vr}{2}} \qquad\mbox{for}\qquad 
-\infty < r < -\infty
$$
form a semigroup under multiplication, i.e.
$$
f_{u_1, v_1}(r)f_{u_2, v_2}(r) = f_{u_1 + u_2, v_1 + v_2}(r)
$$
and the inverse Fourier transform relation \eqn{Nn638}
implies that the measures
$\rho_{u,v}(x) dx$ form a semigroup under convolution, i.e.
$$
P_{u_1, v_1}(x)dx \ast P_{u_2, v_2}(x)dx = P_{u_1 + u_2, v_1 + v_2}(x)dx,
$$
as required. ~~~$\bsq$

This semigroup is {\em holomorphic} in the sense that 
the density functions
$\rho_{u, v}(x)$ are holomorphic  functions of $s = - \frac{u+v}{2}+ ix$
in the cone where they are defined.
Below we compute the moments of the distributions $\rho_{u, v}(x) dx$ using
Lemma~\ref{lem48}.

Before doing that, we derive a formula for the
logarithmic derivatives of the theta function $\theta (t)$ at $t =1$.
Set
\begin{equation}~\label{719a} 
R_k := \frac{d^k}{dt^k} \log\theta(t)|_{t=1}.
\end{equation}

\begin{theorem}~\label{th75}
For each $k \geq 1$, 
\begin{equation}\label{eq751}
R_k \in \QQ[\psi_2],
\end{equation}
where
$$\psi_2 := \pi\theta(1)^4 = \frac{\pi^2}{\Gamma(3/4)^4}.$$
It is given by an even polynomial of degree at most $k$
in $\psi_2$,
\end{theorem}

\noindent\paragraph{Remark.}  In particular, we have
\begin{align}
R_1 &= -\frac{1}{4}\nonumber\\
R_2 &= \frac{1}{8} + \frac{1}{32} \psi_2^2 \nonumber\\
R_3 &= -\frac{1}{8} - \frac{3}{32} \psi_2^2 \nonumber\\
R_4 &= \frac{3}{16} + \frac{9}{32} \psi_2^2-\frac{1}{256} \psi_2^4.\nonumber
\end{align}

\noindent\paragraph{Proof.}
Consider the polynomial ring $\QQ[x(t), y(t), z(t)]$
generated by the three functions 
\begin{align}
x(t)&:=\pi \theta(t)^4          \\
y(t)&:=  -4 \frac{\eta^{(1)}(t)} {\eta (t)} =    
\frac{\pi}{3} E_2(t) \nonumber\\
    &= \frac{\pi}{3} (1 - 24\sum_{1\le m} 
\frac{m e^{-2\pi mt}}{1-e^{-2\pi mt}}).\\
z(t)&:=y(t)+ 4\frac{\theta^{(1)}(t)}{\theta(t)}.
\end{align}
These functions are all (essentially) modular forms of weight $2$ on
the theta group. (The Eisenstein series $E_2(t)$  is not quite a modular
form.)
One has
\[
x(1) = \psi_2, ~~~~~ y(1) = 1, ~~~~~ z(1) = 0,
\]
where the last two are derived using the transformation laws for
$\theta(t)$ and $\eta(t)$. It follows that if $g(t) \in \QQ[x(t), y(t), z(t)]$,
then $g(1) \in \QQ[\psi_2].$

We claim that the polynomial ring $\QQ[x(t), y(t), z(t)]$ is closed
under differentiation $\frac{d}{dt}.$ Indeed one has
\begin{align}
\frac{d}{dt}x(t) &= x(t)(z(t) - y(t)) \\ 
\frac{d}{dt}y(t) &= -\frac{1}{2}y(t)^2 + \frac{1}{24}x(t)^2
   +\frac{1}{8}z(t)^2 \\
\frac{d}{dt}z(t) &= \frac{1}{6} x(t)^2 - y(t) z(t)-
\frac{1}{2} z(t)^2.
\end{align}
These can be deduced using properties of derivatives of modular forms;
the
operator $\frac{d}{dt}+y(t)$ must take both $x(t)$ and $z(t)$ 
into a modular form of weight 4, the specific one being determined
by the first few Fourier coefficients, and $\frac{d}{dt}+\frac{1}{2}y(t)$
does the same for $z(t).$ The procedure is that used in 
Lang~\cite[Chap. 10, Thm. 5.3, p. 161]{La76}.
Furthermore, by symmetry, we observe that the subring
$\QQ[x(t)^2,y(t),z(t)]$ is also closed under differentiation.

Now, we observe that
\[
\frac{d}{dt} \log\theta(t) = \frac{\frac{d}{dt} x(t)}{4x(t)}
=
\frac{1}{4}(z(t)-y(t))\in \QQ[x(t)^2,y(t),z(t)],
\]
so for each $k\ge 1$, $\frac{d^k}{dt^k} \log\theta(t)\in
\QQ[x(t)^2,y(t),z(t)]$.  The theorem follows upon evaluation at $t=1$.
The explicit formulas were found by computation.
~~~$\bsq$

We recall that the {\em cumulants} $\kappa_k= \kappa_k(P)$ 
of a probability distribution $P(dx)$ are
defined in terms of the characteristic function of $P$ by:
\begin{equation}~\label{722a}
\kappa_k := i^{-k} \left(\frac{d^k}{dz^k} \log(\int_{-\infty}^\infty e^{iz x}
P(dx))\right)_{z=0}.
\end{equation}
In particular, $\kappa_1(P)$ is the mean and $\kappa_2(P)$ is the variance of
$P$.  We extend this definition to the complex measures $\rho_{u,v}(x) dx$.
The usefulness of cumulants 
for infinitely divisible distributions (as opposed to moments)
is that they scale nicely with the parameter $u$.

\begin{theorem}~\label{th76}
For real $u >0$ and real $|v| < u$,
the mean value 
$$ \kappa_1(\rho_{u,v}) = -\frac{iv}{2}. $$
For all  $k \geq 2$,
the $k$-th cumulant of the distribution $\rho_{u,v}(x)dx$ has the form
$$\kappa_k(\rho_{u,v}) = c_k u,$$
where $c_k$ in $\QQ[\psi_2]$,
with
$$\psi_2 := \pi\theta(1)^4 = \frac{\pi^2}{\Gamma(3/4)^4}.$$
Moreover, if $k \ge 3$ is odd, then $c_k=0$, while if $k \geq 2$
is even, then $c_k$ is an even polynomial of degree $k$ in $\psi_2.$
\end{theorem}

\noindent\paragraph{Remark.}  In particular, we have
\begin{align}
c_2 &=  -\frac{1}{2} + \frac{1}{8}\psi_2^2 = 4(R_1 + R_2)\nonumber\\
c_4 &=  -1 + \psi_2^2 + \frac{1}{16}\psi_2^4. \nonumber
\end{align}

\noindent\paragraph{Proof.}
We have:
\begin{align}
\int_{-\infty}^\infty e^{iz x} \rho_{u,v}(x) dx
&= \frac{1}{2\pi}
 \theta(1)^u 
\int_{-\infty}^\infty e^{iz x} Z_\QQ(-u,-\frac{u+v}{2}+ix) dx\\
&= \frac{1}{2\pi}
\theta(1)^u 
e^{(u+v)z/2}
\int_{-\infty}^\infty e^{z(-\frac{u+v}{2}+ix)} Z_\QQ(-u,-\frac{u+v}{2}+ix) dx.
\end{align}
Applying Lemma~\ref{lem48} with $Q(T)=1$, $\cw=-u$, and
$\sigma=-\frac{u+v}{2}$, replacing $z$ by $-z$,
and multiplying by $\theta(1)^ue^{(u+v)z/2}$,
we obtain
\[
\int_{-\infty}^\infty e^{iz x} \rho_{u,v}(x) dx
=
\theta(1)^u
e^{(u+v)z/2}
\theta(e^{-2z})^{-u},
\]
and thus
\[
\log(\int_{-\infty}^\infty e^{iz x} \rho_{u,v}(x) dx)
=
\frac{vz}{2} + u (\frac{z}{2}+\log\theta(1)-\log\theta(e^{-2z})).
\]
We expand this in a Taylor expansion in $z$; clearly, $v$ contributes only to
the first cumulant.  Conversely, since $u$ is multiplied by an even
function of $z$, it contributes only to the even cumulants; the claim
for $\kappa_1$ follows immediately.
Since
\[
\log\theta(1+x) = \log\theta(1) + \sum_{k\ge 1} R_k x^k/k!,
\]
we find that
\[
\log\theta(1)-\log\theta(e^{-2z}) = -\sum_{k\ge 1} R_k (e^{-2z}-1)^k/k!,
\]
and thus has Taylor coefficients in $\QQ[\psi_2]$.
~~~$\bsq$

Since the moments of a measure are polynomials in its cumulants, we obtain:

\begin{cor}
For real $u>0$ and $|v|<u$ and integer $k\ge 0$, the moments
\[
M_k(u,v)=\int_{-\infty}^\infty x^k \rho_{u,v}(x) dx
\]
satisfy $M_k(u,v)\in \QQ[\psi_2][u,v]$, where
$\psi_2 := \pi\theta(1)^4 = \frac{\pi^2}{\Gamma(3/4)^4}.$
\end{cor}

We determine the mass of the Feller
canonical measure.

\begin{theorem}~\label{thm79}
For all $\cw \in \CC$,
\[
\frac{1}{2\pi} \int_{-\infty}^\infty \xi_\QQ(\cw,\frac{\cw}{2}+ix) dx
=
\theta(1)^\cw (\frac{\psi_2^2}{16}   -\frac{\cw}{8} -\frac{1}{4}).
\]
In particular, the mass of the Feller canonical measure $M\{dx\}$is
\begin{equation}~\label{729}
\int_{-\infty}^\infty M(x) dx = \frac{\psi_2^2}{8} - \frac{1}{2},
\end{equation}
which is $\frac{\pi^2}{8} \theta(1)^8 -\frac{1}{2} \approx 1.8946.$
\end{theorem}

\noindent\paragraph{Proof.}
We deduce this using Lemma~\ref{lem48}.  From the uniform Schwartz
property of $\xi_\QQ$, we see that the left-hand side is entire in $\cw$.
It thus suffices to prove the theorem when $\Re(\cw)<0$.  We then have
\begin{align}
\frac{1}{2\pi} \int_{-\infty}^\infty \xi_\QQ(\cw,\frac{\cw}{2}+ix) dx
&=
\frac{1}{2\pi}
\int_{-\infty}^\infty Q(\frac{\cw}{2}+ix) Z_\QQ(\cw,\frac{\cw}{2}+ix) dx\\
&=
Q(-\frac{d}{dz})(\theta(e^{-2z})^{\cw})|_{z=0},
\end{align}
where $Q(z) = z(z -\cw)/2\cw$.  Differentiating and evaluating, we obtain:
\[
\frac{1}{2\pi} \int_{-\infty}^\infty \xi_\QQ(\cw,\frac{\cw}{2}+ix) dx
=
\theta(1)^{\cw}
\left((2R^2 +(\cw-2)R_1 +2\cw R_1^2 \right)
=
\theta(1)^w (\frac{\psi_2^2}{16} -\frac{w}{8} - \frac{1}{4}),
\]
as required.
Finally, since  $M(x) = \frac{1}{\pi} \xi_{\QQ}(0, -ix),$ we obtain the
mass of the Feller canonical measure on taking $\cw = 0.$
~~~$\bsq$

%
%
%

\subsection{Location of Zeros}

We  know little about the
location of the zeros of 
$\xi_{\QQ}(-u, s)$, for negative real $-u$, aside from the
general bound given by Lemma~\ref{le41}.
We raise the following questions.

\noindent{\em Question 1.} Is $\Re(Z(\cw, s)) > 0$ in
the entire real-codimension one cone 
$\sC^{-} := \{ (\cw, s) : \cw = u \in \RR,~ u < \Re(s) < 0  \}$?

If true, this would extend the result of Theorem~\ref{nth61}
to exclude zeros from the open cone. At $u=0$ all the
zeros are strictly outside the cone, so to prove this 
result it would suffice to show that there are never any
zeros on the boundary of the cone.

\noindent{\em Question 2.} For each fixed $u < 0$ are the zeros
of $\xi_{\QQ}(u, s)$ confined to a vertical strip $|\Re(s)| < g(u)$
for some function $g(u)$?

The result of \S5 shows that this is true on the boundary
plane $u =0;$ perhaps it persists in the region  $u < 0.$ 
The results of \S5 also suggest for $u >0$ a limited movement of zeros
in the vertical direction as $u$ varies, so we 
ask the following question for negative $u$.

\noindent{\em Question 3.} For fixed $u < 0$ let
$N_{u}(T)$ count (with multiplicities)
 the total number of zeros $\rho = \beta + i \gamma$ of
$\xi_{\QQ}(u, s)$ lying in the horizontal strip  $|\Im(\rho)| \leq T.$ 
Is $N_{u}(T)$
finite for each $T >0$, and if so, does it obey the same 
asymptotic formula as that for $u > 0$ in Theorem~\ref{SPth51}?
Here we only ask for a remainder term smaller than the main
term, possibly with a different dependence on $u$.

For real $-u < 0$, Theorem~\ref{nth61} shows that
there are no zeros on  the center line $\Re(s) = -\frac{u}{2}.$
Thus 
the nonreal zeros always occur in quadruples $\{ \rho, -\frac{u}{2} - \rho,
\bar{\rho}, -\frac{u}{2}- \bar{\rho}\}.$

\noindent{\em Question 4.} The Riemann hypothesis is encoded in the
location of the zeros of $\xi_{\QQ}(0, s),$ asserting they are on the four
lines $\Re(s) = \pm 1,\pm 2.$ Indeed, the assertion that $\xi_{\QQ}(0,s)$
has no zero with $|\Re(s)|<1$ is equivalent to the Riemann hypothesis.
Does the convolution semigroup structure for $u < 0$ play any factor in
controlling the location of these zeros?

This question can be studied by considering more general
convolution semigroups.

%
%
%

\section{Location of Zeros: Complex $\cw$}
\hsp
The zero locus $\sZ_{\QQ}$ of the entire function $\xiq(\cw,s)$ 
decomposes into a countable union of Riemann
surfaces embedded in $\CC^2$; we call these {\em components}.
How many components are there in $\sZ_{\QQ}$?
While we cannot answer this question, we can at least show
that certain zeta zeros (on the slice $\cw = 1$)
lie on the same component.

\begin{lemma}\label{Nle71}
The Riemann zeta zeros $\rho_7 \approxeq \frac{1}{2} + 42.04i$ and 
$\rho_8 \approxeq \frac{1}{2} + 42.90i$ appearing
as zeros of $\zq (\cw, s)$
in the slice $\cw =1$ belong 
to the same component of $\sZ_{\QQ}$.
\end{lemma}

\paragraph{Proof.}
We deform $\cw= u$ through real values in the interval $1 \le u \le 2$.
By numerical computation we
find that at $u= \frac{3}{2}$ these two zeros have moved off the critical 
line to assume complex conjugate values.

We use the symmetry that when
 $\cw = u$ is real, if $\rho$ is a zero, then so are $\bar{\rho}$, 
$1- \rho$ and $1- \bar{\rho}$.
In particular, zeros cannot move off the critical line except 
by combining in pairs.
As $u$ changes, at some point $u_0$ they must coalesce on the 
critical line as a double zero, then as $u$ changes go off the line, becoming
a pair of complex conjugate zeros.
The point of coalescence at $u_0$ of two zeros could be either the 
intersection of two
different components of $\sZ_{\QQ}$ (the intersection having 
real codimension 4 in $\CC^2$) or a single 
component of $\sZ_{\QQ}$ having a branch point of order two there 
(when viewed as projected
on the $\cw$- plane.)
The latter case must occur, because in the first case the movement of the zeros
$\rho (u)$ would have a first derivative as a function of $u$ which
varies analytically in $u$ at the critical point.
This manifestly does not happen, because as a function
of $u$ the zeros first move vertically
on the critical line, then change directions at $u_0$ to move horizontally
off the line.
Thus the component forms a single Riemann surface, 
with a path on it from
$(u,s) = (1, \rho_7 )$ to $(1, \rho_8)$.~~~$\bsq$

It seems reasonable to guess that the zeta zeros 
$\{ \rho = \sigma + it : {\rm Im} (t) > 0 \}$ lie on the same component of 
$\sZ_{\QQ}$.
If so, the same would hold for $\{ \rho = \sigma + it : {\rm Im} (t) < 0 \}$,
since the zero set is invariant under complex conjugation,
i.e. $\bar{\sZ}_{\QQ} = \sZ_{\QQ}.$
The simplest hypothesis concerning the zero set would
seem to be that it is the closure of a single irreducible complex-analytic
variety of multiplicity one. However we do not have any strong evidence for
this hypothesis.

%
%
%

\section{General Number Fields}

We now briefly consider the Arakelov zeta function
for a general algebraic number field $K$. Many of the results
extend to general $K$ but the positivity property of
Theorem~\ref{th11} does not.

In the case of 
the Gaussian field $\QQ(i)$, 
we have the identity
\begin{equation}~\label{109a}
Z_{\QQ(i)}(\cw,s) = 2 \zq(2\cw, 2s),
\end{equation}
derived in the appendix. Thus all the results proved here immediately apply 
to $K = \QQ(i).$

For a general algebraic number field $K$, the  
Arakelov two-variable zeta function $\zk(\cw, s)$ 
has a functional equation. Furthermore it can be shown that
$\xi_K(\cw, s) := \frac{s(s -\cw)}{2\cw} \zk(\cw, s)$
is an entire function of order one and infinite type of two variables,
by generalization of the proofs for $K=\QQ$.

The  proofs given here for  the distribution of zeros of $\zq(\cw, s)$
for positive real $\cw$ partially extend to general $K$. 
The proofs of confinement of the zeros to a vertial strip of
width depending on $\cw$ extend to a few fields such as 
the Gaussian field $\QQ(i)$; they depend on the existence of
an associated
 Dirichlet series with a nonempty half-plane of absolute convergence.
For general $K$ the Dirichlet series has
a nonempty half-plane of convergence for $\cw$ a positive integer,
but perhaps not for any other values of $\cw$.
One expects to get estimates for zeros to height $T$ for such
integer values of $\cw$ .
We do not know whether for fixed positive noninteger
real $\cw$ and general $K$ the zeros are confined to a vertical strip
of finite width,
or that a generalization of the zero counting bounds in
Theorem~\ref{SPth51} holds, counting zeros in
a horizontal strip $\Im(z) < T$.
The latter seems plausible, because the zeros are confined at 
positive integer $\cw$, but if so, new proofs are needed.

The convolution semigroup property, 
of a family
of complex-valued measures on a cone associated to negative real $\cw$, 
continues to hold for general
number fields $K$. 
The associated measures are real-valued
on the ``critical line.''
However, the  proof of positivity of such measures 
on the ``critical line''
for $\QQ$ given in  Theorem~\ref{th11}
extends only to a few specific number fields, such
as $\QQ(i).$ Our proof used a product
formula for the associated modular form, which permits
calculations with its logarithm and yields an explicit
form for the associated Feller canonical measure. 
Such product forms exist only for modular forms 
all of whose zeros are at cusps.
The modular forms associated to most imaginary quadratic fields
generally do not have a product formula, because the associated
modular forms have zeros in the interior of a fundamental domain,
and the logarithm of such forms are multivalued functions. 

There are imaginary quadratic number fields 
with class number one for  which the positivity property
does not hold.
One can show that the  positivity property holds for a field $K$
if and only if 
$\xi_K(0, it)$ is nonnegative for all real $t$, where
$\xi_K(\cw, s)= \frac{s(s- \cw)}{2\cw}Z_K(\cw, s).$
For $K= \QQ(\sqrt{-2})$ one finds that
$$\xi_K(0, it) = 2\cos (t\log \sqrt{2})\xi_{\QQ}(0, 2it), $$
which clearly has sign changes.
We also found by  computer calculation that 
$\xi_K(0, it)$ for $K = \QQ(\sqrt{-11})$ changes sign
between $t=3.10$ and $3.15$, and in addition on the line $\cw=-1$
there is a sign change, with $Z_K(-1, 4i) < 0.$
For $K=\QQ(\sqrt{-19})$ there is a  sign change of $\xi_K(0, it)$
between $t=2.0$ and $t=2.1$. 
It remains conceivable that there exist 
imaginary quadratic fields having  the
positivity property,  whose
associated modular form does not have a product formula.
In support of this,
computer experiments for the imaginary quadratic
number fields $\QQ(\sqrt{-3})$ and $\QQ(\sqrt{-7})$
did not locate any sign changes for $\xi_K(0, it)$.

In a different direction, the positivity property of the convolution semigroup
for negative real $u$ on the ``critical line'' generalizes  
to certain classes 
of modular forms not associated to number fields,
which do have a product formula,
as we hope to treat elsewhere.

%
%
%

\newpage
\appendix
\section{Appendix: Arakelov zeta function of van der Geer and Schoof}

This appendix summarizes the framework of van der Geer
and Schoof~\cite{vdGS99}, and obtains explicit formulas for
the two-variable zeta function for $K= \QQ$ and $\QQ(i).$
The expression of van der Geer and Schoof for the
Arakelov two-variable  zeta function  is, formally,
\begin{equation}\label{A101}
\hzk (\cw, s) \cong \int_{Pic(K)} e^{s h_0 ([D]) + (\cw -s) h_1 ([D])} d[D] ~.
\end{equation}
In this expression $h^0 ([D])$ resp. $h^1 ([D])$ are analogous to 
the ``dimension'' of a sheaf cohomology group.
They give a direct definition of $h^0 ([D])$, 
and then define $h^1 ([D])$ 
indirectly~{\footnote
{A. Borisov~\cite{Bo98} has given a direct
definition of $h^1 ([D])$ in some cases, with a proof of the
formula \eqn{A102}.}}
 to be 
\begin{equation}\label{A102}
h^1 ([D]) := h^0 ( [\kappa_K] - [D]) ~,
\end{equation}
where $\kappa_K$ is the ``canonical'' Arakelov divisor 
for the ring of integers of $K$, which is
what the Riemann-Roch formula predicts. 
We now define $Pic(K)$ and $h^0([D]).$
The value $h^0 ([D])$
turns out to be the logarithm of a multivariable theta function at a
specific point depending on $[D]$, see \eqn{A113}.

In what follows, let
$K$ be an algebraic number field, with  $O_K$ its ring of integers
and $\Delta_K$ its discriminant.
Set  $[K : \QQ] = n = r_1 + 2r_2$, with $r_1$ real places and
$r_2$ complex places. We denote archimedean places of $K$ by $\sigma$
and non-archimedean places by $\nu.$ 

\begin{defi}~\label{deA11}
{\em
(i) An {\em Arakelov divisor} $D$ is a formal finite sum over
the non-archimedean places $\nu$ of $K$ and the 
$r_1 + r_2$ archimedean places $\sigma$,
$$ D = \sum_{\nu < \infty} n_{\nu} \nu ~~+ 
\sum_{\sigma} x_\sigma \sigma $$
in which each $n_\nu$ is an integer and each $x_\nu$ is a real number
(even at a complex place $\sigma$.)

(ii) An Arakelov divisor $D$ is {\em principal} if there is an
element $\alpha \in K^*$ such that
$$D = (\alpha) = \sum_{\nu < \infty} ord_{\nu}(\alpha) \nu
+ \sum_{\sigma} x_\sigma(\alpha) \sigma, $$
in which $x_\sigma(\alpha)$ equals $ \log |\sigma(\alpha)|$ 
or $2 \log |\sigma(\alpha)|$ according as $\sigma$ is a real place
or a complex place. Here $\sigma(.)$ runs over all embeddings of
$K$ into $\CC$, with the convention that only one out of 
each  conjugate complex pair of complex embeddings is used. 

(iii) $Div(K)$ denotes the
 group of Arakelov divisors (under addition). 
The {\em Arakelov divisor class group} $Pic(K)$ is the
quotient group by the subgroup of principal Arakelov divisors.
The divisor class of $D$ is denoted $[D]$.
}
\end{defi}
The roots of unity $\mu_K$  in $K$ have Arakelov 
divisor zero. They fit in the exact sequence
\begin{equation}~\label{A103}
0 \to \mu_K \to K^* \to Div(K) \to Pic(K) \to 0.
\end{equation}

\begin{defi}~\label{deA2}
{\em 
The {\em degree} $deg(D)$ of an Arakelov divisor is the
real number
$$ deg(D) : = \sum_{\nu < \infty} n_{\nu} \log N\nu + 
\sum_{\sigma } x_{\sigma}. $$
Here $N \nu := |O_K /P_{\nu}|$, where 
$P_\nu = \{ \alpha  \in O_K : |\alpha|_{\nu} < 1\}$,
in which $N\nu$ is the number of elements in the residue
field of $\nu$.
}
\end{defi}

Principal divisors have degree zero, so the degree $\deg([D])$ 
is well-defined on Arakelov divisor classes.

\begin{defi}~\label{deA3}
{\em 
The {\em norm} $N(D)$ of an Arakelov divisor $D$ is
$$N(D) = \exp (deg(D)) = \prod_{\nu < \infty} (N\nu)^{n_\nu}
\prod_{\sigma} e^{x_\sigma}. $$

}
\end{defi}

\begin{defi}~\label{deA4}
{\em 
(i) The {\em ideal} $I_D$ associated to  an Arakelov divisor $D$ 
at the finite places is the fractional ideal
$$I_D := \prod_{\nu < \infty} P_\nu^{- n_\nu}, $$
where $P_{\nu}$ denotes  the prime ideal at $\nu$.

(ii) The lattice structure 
associated to an  Arakelov divisor $D$
at the infinite places is a positive
 inner product on $\RR^{r_1} \times \CC^{r_2}$.
defined as follows.
At a real place, $x_\sigma$ determines a scalar product on $\RR$
such that $||1||^2 = e^{-2x_\sigma}.$ At a complex place
$x_\sigma$ determines a Hermitian inner product on $\CC$
such that $||1||^2 = 2e^{-x_\sigma}.$ 
The combined inner product is 
$$||(z_\sigma)||_D^2 := \sum_{\sigma} |z_\sigma|^2 ||1||_\sigma^2. $$
The {\em (metrized) lattice} $\Lambda_{D}$ associated to $D$ 
is the fractional ideal $I_{D}$ (viewed as a subset of $K$) embedded
in $\RR^{r_1} \times \CC^{r_2}$ as Galois
conjugates of each element $\alpha$, with distance function
measured by this inner product.
}
\end{defi}

The number field $K$ embeds as a dense subset of $\RR^{r_1} \times \CC^{r_2}$,
while each fractional ideal $I_D$ embeds discretely as a lattice in
this space.
The archimedean coordinates $x_\sigma$ define a metric structure 
at the infinite places such that 
$$ Covol(\Lambda_D) = \frac{\sqrt{\Delta_K}}{N(D)}, $$
where $\Delta_K$ is the discriminant of $K$. The 
Arakelov class group $Pic(K)$ parametrizes isometry classes of
lattices that have compatible $O_K$-structures under multiplication.
Following Szpiro, the  {\em Euler-Poincar\'{e} characteristic} $\chi(D)$ of an
Arakelov divisor $D$ is defined as 
$$
\chi(D) := - \log( covol(\Lambda_{D})) = deg(D) - \frac{1}{2} \log \Delta_K. 
$$
It is well-defined on divisor classes $[D]$.
In general the Arakelov class group
\begin{equation}~\label{A105a}
Pic(K) \simeq Cl(K) \times \RR \times \TT^{r_1 + r_2 -1}, 
\end{equation}
where $Cl(K)$ denotes the (wide) ideal class group of $K$,
and  the second factors combined are  $\RR^{r_1 + r_2}/U(K)$
where $U(K)$ is an $r_1 + r_2 -1 $ dimensional lattice
given by logarithms of (absolute values of) Galois conjugates of units.
Note that  $Pic^0(K)$, the group of divisor classes of
degree $0$, is compact, and its volume is $h_K R_K$ where
$h_K$ and $R_K$ are the class number and regulator of $K$,
respectively.

\begin{defi}~\label{deA5}
{\em 
The {\em canonical divisor} $\kappa_K$ of a number field $K$ 
is the Arakelov
divisor whose associated ideal part $I_{\kappa_K}$ is the inverse different
$\sgd_K^{-1}$ for $K/\QQ$, 
and all of whose archimedean coordinates $x_\sigma = 0$.
}
\end{defi}

These definitions imply that 
$$
N(\kappa_K) = N(I_{\kappa_K})^{-1} = N(\sgd_K) = \Delta_K.
$$

\begin{defi}~\label{deA6}
{\em
(i) An Arakelov divisor $D$ is {\em effective} if $O_K \subseteq I_D.$

(ii) The {\em effectivity} $e(D)$ of an  effective divisor $D$ is
$$ e(D) = \exp ( - \pi ||1||_D^2), $$
in which
$$ ||1||_D^2 := \sum_{\sigma~~ real} e^{-2 x_\sigma} + 
\sum_{\sigma~~ complex} 2 e^{- x_\sigma}.$$
The {\em effectivity} $e(D)$ of a non-effective divisor is $0$.
}
\end{defi}

The effective divisors are those Arakelov divisors
with each $n_\nu \geq  0$ and the effectivity $e(D)$ of any
divisor takes
a value $0 \leq e(D) < 1.$  
The only ``functions''
$ \alpha \in K^*$ whose associated  principal divisors $(\alpha)$ 
are effective
are the roots of unity in $K$, whose associated Arakelov divisor is $0$,
the identity element. By convention we add a  symbol $(0)$ to  represent
a ``divisor at infinity'' corresponding to the element $0 \in K$,
with the convention that $ (0) + D = (0)$ for all Arakelov divisors $[D]$
and we define the  effectivity $e( (0)) := 1.$

\begin{defi}~\label{deA7}
{\em
(i) The {\em order } $H^{0}(D)$ of the group of effective
divisors associated to an Arakelov
divisor  $D$ is
$$ H^0(D) := \sum_{\alpha \in I_D} e( (\alpha) + D). $$
This sum includes a term  $\alpha = 0 \in K$, with the convention that
$e( (0) + D)) := 1,$ so that $H^0([D]) \geq 1.$

(ii) The {\em effectivity dimension} $h^{0}(D)$ of $H^0(D)$ is given by
$$ h^0(D) = \log H^0(D). $$
One has $h^0(D) \geq 0$. 
}
\end{defi}

The quantities $H^0(D)$ and $h^{0}(D)$ are constant for all
divisors in an Arakelov divisor class $[D] \in Pic(K)$ and
may therefore be denoted $H^0([D])$ and $h^0([D])$, respectively.
van~der Geer and Schoof \cite[Prop. 1]{vdGS99} state the following
result.

\begin{theorem}~\label{thA1}
{\bf (Riemann-Roch Theorem for Number Fields)}
For any algebraic number field $K$ and any Arakelov
divisor class $[D] \in Pic(K)$, 
\begin{equation}\label{A104a}
h^{0}([D]) - h^{0}([\kappa_K] - [D]) = 
\deg([D]) - \frac{1}{2} \deg([\kappa_K]),
\end{equation}
in which $\kappa_K$ is the canonical Arakelov divisor for $K$,
and $\deg([\kappa_K])= \log \Delta_K$.
\end{theorem}

van der Geer and Schoof
defined a new invariant $\eta(K)$ of a number field $K$, 
and an Arakelov analogue of the  genus $g(K)$ of $K$.

\begin{defi}~\label{deA8}
{\em
(i) The invariant $\eta(K)$ of $K$ is defined by
\begin{equation}~\label{A105b}
 \eta(K) := H^0(O_K).
\end{equation}

(ii) The {\em genus} $g(K)$ of $K$ is defined by
\begin{equation}~\label{A106} 
g(K) := h^0(\kappa_K) = h^0(O_K) + \frac{1}{2} \log \Delta_K
 = \log (\eta(K) \sqrt{\Delta_K}).
\end{equation}
}
\end{defi}

For the rational number field
\begin{equation}~\label{A107} 
\eta(\QQ) = \omega := \frac{\pi^{1/4}}{\Gamma(\frac{3}{4})} \approx 1.08643,
\end{equation}
and 
$g(\QQ) = \log \eta(\QQ) \approx  0.0829015.$
The value $\omega = \theta(1)$, where
$\theta(\cdot)$ is the theta function \eqn{S104}.
One also has
$$
\eta(\QQ(i)) = \big(\frac{2 + \sqrt{2}}{4} \big)\omega^2 ,
$$
see \cite[pp. 16-17]{vdGS99}.

The genus of a function field is usually defined to be
the dimension $h^0(\kappa)$ of the vector space of
effective divisors for the canonical class
$\kappa$; this motivates the  definition of 
the genus $g(K)$ of a number field.
For a function field the degree $deg(\kappa)$ 
of the canonical class is $2g-2$,
so one may consider 
$$\tilde{g}(K) := \frac{1}{2}(\deg(\kappa_K) + 2)$$
as a second analogue of genus for a number field. 
This analogue appears on the right side of
 the Riemann-Roch theorem for number fields.  One has
\begin{equation}~\label{A107a}
\tilde{g}(K) = 1 + \frac{1}{2} \log \Delta_K.
\end{equation}
In particular  $\tilde{g}(\QQ)=1$, and
$\tilde{g}(\QQ(i)) = 1 + 2 \log 2$.
The two notions of genus agree in the
function field case and differ in the number field case.

Below we obtain explicit integral
formulas for the two-variable  zeta functions for  $K= \QQ$
and $\QQ(i)$, of the form \eqn{A101}, and  also indicate the
form of the two-variable zeta function for a general algebraic
number field $K$.  
 
To define an integral over the Arakelov class group, one must
specify a measure on the group. 
For compact groups it is Haar measure,
and on noncompact additive groups $\RR$ it is $dx.$
It is convenient to replace additive groups by multiplicative
group $\RR_{> 0}$ using the change of variable
$y_\sigma = e^{- x_\sigma}$ at real places and the
appropriate measure becomes $\frac{dy}{y}.$
For $\QQ$ and $\QQ(i)$ the Arakelov class group is 
isomorphic to $\RR.$

\noindent\paragraph{Case $K= \QQ$.}
There is a single real place $\sigma$. For
representatives of Arakelov divisor classes $[D]$ we may
take $D$ to have ideal $O_{\QQ} = \ZZ$ and with value
at the infinite
place $x_\sigma \in \RR$ arbitrary, with measure $dx$
at the infinite place.  Thus the Arakelov
class group $Pic(\QQ) \equiv \RR$, the additive group,
with $x \in \RR$ being the degree of the divisor.
We have
\begin{eqnarray}~\label{A108}  
H^0([D]) & = & \sum_{n \in \ZZ} e( (n) + D)     
\nonumber \\
& = &  1 + 
\sum_{n \in \ZZ \backslash \{0\}}  \exp (- \pi e^{2\log |n| - 2x_\sigma})
\nonumber \\
& = &   \sum_{n \in \ZZ}    e^{- \pi n^2 y_{\sigma}^2},
\end{eqnarray}
in which one sets $y_\sigma = e^{- x_\sigma}$. 
Using
the multiplicative change of variable $y_\sigma = e^{- x_\sigma}$
we  identify  $Pic(\QQ)$ with the multiplicative group 
$ \equiv \RR_{> 0}$
with measure $\frac{dy}{y}.$ Thus we have
$$H^0 ([D]) = \theta( y_\sigma^2), $$
in terms of the theta function \eqn{S104}.

The different $\sgd_{\QQ} = (1),$  so the canonical divisor
$\kappa_{\QQ} = 0.$ Consequently we have
$$H^{1}([D]) := H^0 ( -[D]) = \theta (\frac{1}{y_{\sigma}^2}). $$
We obtain
\begin{eqnarray}~\label{A109}
Z_{\QQ}(\cw,s) & := &  \int_{Pic(\QQ)} e^{s h_0 ([D]) + (\cw-s) h_1 ([D])} d[D]
\nonumber \\
& = & \int_{Pic(\QQ)} H^0([D])^s H^1([D])^{\cw - s} d[D] \nonumber \\
& = & \int_{0}^{\infty} \theta (y^2)^{s} \theta(\frac{1}{y^2})^{\cw - s} 
\frac{dy}{y}.
\end{eqnarray}

\noindent\paragraph{Case $K= \QQ(i)$.}
There is a single complex  place $\sigma$.
The class number of $O_K = \ZZ[i]$ is one, so all Arakelov
divisor classes $[D]$ in $Pic(\QQ(i))$ have a representative 
$D = x_\sigma \sigma,$
whose associated ideal is $O_K$. Thus $Pic(K) \equiv \RR$
as an additive group. Letting $\alpha= m + ni \in \ZZ[i]$,
we have
\begin{eqnarray}~\label{A110}
H^0([D]) & = & \sum_{\alpha \in \ZZ[i]} e( (\alpha) + D) =    
  1 + \sum_{\alpha \in \ZZ[i] \backslash \{0\}}  
\exp(- 2\pi e^{2\log|\alpha| - x_\sigma)})
\nonumber \\
& = &   \sum_{m\in \ZZ} \sum_{n \in \ZZ}   e^{- 2\pi (m^2 +n^2)y_{\sigma}},
\nonumber \\
& = & \left(\sum_{m\in \ZZ} e^{-2 \pi m^2 y_{\sigma}} \right)^2.
\end{eqnarray}
in which $y_\sigma = e^{-x_\sigma}.$ Thus
\begin{equation}~\label{A11}
H^0([D]) = \theta(2y_\sigma)^2.
\end{equation}

The different $\sgd_{\QQ(i)} = ((1 + i)^2) = (2),$ and the
canonical divisor $\kappa_{\QQ(i)} = \nu_{(1+i)}^{-2} = (\frac{1}{2}). $
Consequently we have
$$H^{1}([D]) := H^0 ( \kappa_{\QQ(i)} - [D]) = 
\theta (\frac{1}{2 y_{\sigma}}).$$
We obtain
\begin{eqnarray}~\label{A111}
Z_{\QQ(i)}(\cw,s) & := &  \int_{Pic(\QQ(i))} e^{s h_0 ([D]) + 
(\cw-s) h_1 ([D])} d[D]
\nonumber \\
& = & \int_{Pic(\QQ(i))} H^0([D])^{s} H^1([D])^{\cw - s} d[D] \nonumber \\
& = & \int_{0}^{\infty} \theta (2y)^{2s} \theta(\frac{1}{2y})^{2\cw - 2s} 
\frac{dy}{y} \nonumber \\
& = & 2 \int_{0}^{\infty} \theta (t^2)^{2s} \theta(\frac{1}{t^2})^{2\cw - 2s}
\frac{dt}{t} 
\end{eqnarray}
using the change of variables $t^2 = 2y$. Comparing this with \eqn{A109}
yields
\begin{equation}~\label{A112}
 Z_{\QQ(i)}(\cw,s) = 2Z_{\QQ}(2\cw, 2s). 
\end{equation}

\noindent\paragraph{General Number Fields.}
Let $K$ be an algebraic number field, of degree $[K:\QQ]=n$.
We follow Lang~\cite[Chapter 13]{La94} for Hecke's
functional equation for the Dedekind zeta function.
One can show that
\begin{equation}~\label{A113}
Z_K(\cw,s) = \frac{2}{w(K)} \sum_{\sga \in Cl(K)} \int_{0}^\infty 
\left( \int_{E} \theta(t^{2/n}\sgc, \sga )^s 
\theta(t^{-2/n}\sgc^{-1}, \sgd_K \sga^{-1})^{\cw - s} d^{*}\sgc \right) 
\frac{dt}{t},
\end{equation}
which uses a decomposition of the Arakelov class group \eqn{A105a}.
Here $\sga $ runs over a set of 
representatives of the (wide) ideal class group,
$w(K)$ counts the number of roots of unity in $K$, and $E$ is a fundamental
domain in the (logarithmic) space of units, with Haar measure $d^{*}\sgc$.
The theta function
$\theta(t^{2/n}\sgc, \sga)$
is defined in Lang~\cite[p.253]{La94} and 
satisfies the functional equation 
\begin{equation}~\label{A115}
\theta(t^{2/n}\sgc, \sga)= \frac{1}{t} 
\theta(t^{-2/n}\sgc^{-1}, \sgd_K \sga^{-1}),
\end{equation}
using the fact that $||\sgc||=1,$ see Lang~\cite[p.257]{La94}.
Using this functional equation and the substitution $u=t^{-2}$  one obtains
\begin{equation}~\label{A116}
Z_K(\cw,s) = \frac{1}{w(K)} \sum_{\sga \in Cl(K)} \int_{0}^\infty 
\left( \int_{E} \theta(u^{1/n}\sgc^{-1},\sgd_K \sga^{-1})^w d^{*}\sgc \right)
u^{s/2} \frac{du}{u}.
\end{equation}
One has the functional equation
$$
Z_K(\cw,s)= Z_K(\cw, \cw-s).$$
For $\cw=1$ one recovers the completed Dedekind zeta function
\begin{equation}~\label{A117}
Z_K(1,s)= \hat{\zeta}_K(s) := 
A_K^{s}\Gamma(\frac{s}{2})^{r_1} \Gamma(s)^{r_2}\zeta_K(s),
\end{equation}
in which $A_K := 2^{-r_2} \pi^{-n/2} \Delta_K^{1/2}   $, 
see Lang~\cite[p. 254]{La94}. \\


\noindent{\tt email:\\
\begin{tabular}{ll}
jcl@research.att.com \\
rains@research.att.com \\
\end{tabular}
 }

\end{document}